\documentclass{amsart}

\usepackage{amsmath}
\usepackage{amscd}
\usepackage{amssymb} 

\newcommand{\cal}{\mathcal}
\newcommand{\bk}{{\bf k}}
\newcommand{\bs}{{\bf s}}
\newcommand{\fg}{{\frak g}}
\newcommand{\fn}{{\frak n}}
\newcommand{\ft}{{\frak t}}
\newcommand{\wh}{\wedge_h}
\newcommand{\ww}{\wedge_{w}}

\DeclareMathOperator{\End}{End}
\DeclareMathOperator{\Img}{Im}

\DeclareMathOperator{\Ker}{Ker}

\newtheorem{theorem}{Theorem}[section]
\newtheorem{theorem/definition}{Theorem/Definition}[section]
\newtheorem{Theorem}{Theorem}
\newtheorem{proposition}{Proposition}[section]
\newtheorem{lemma}{Lemma}[section]

\newenvironment{remark}{\medskip 
\noindent {\bf Remark.}}{\mbox{}}

 \newenvironment{example}{\medskip 
\noindent {\bf Example.}}{\mbox{}} 
\newenvironment{definition}{\medskip 
\noindent {\bf Definition.}}{\mbox{}} 

\begin{document}
\title
{On Quantum de Rham Cohomology}
\author{Huai-Dong Cao \& Jian Zhou}
\address{Department of Mathematics\\
Texas A \& M University\\
College Station, TX 77843}
\email{cao@math.tamu.edu,  zhou@math.tamu.edu}
\begin{abstract}
We define quantum exterior product $\wedge_h$ 
and quantum exterior differential $d_h$
on  Poisson manifolds,
of which symplectic manifolds are an important class of
examples. 
Quantum de Rham cohomology is defined as the cohomology of $d_h$.
We also define quantum Dolbeault cohomology.
Quantum hard Lefschetz theorem is proved.
We also define a version of quantum integral,
and prove the quantum Stokes theorem. 
By the trick of replacing $d$ by $d_h$ and 
$\wedge$ by $\wh$ in the usual definitions, 
we define many quantum analogues of
important objects in differential geometry, e.g. 
quantum curvature. 
The quantum characteristic classes are then studied 
along the lines of classical Chern-Weil theory,
i.e., they can be represented by expressions of quantum curvature.
Quantum equivariant de Rham cohomology is defined 
in a similar fashion.
Calculations are done for some examples, which show that
quantum de Rham cohomology is different from the quantum cohomology
defined using pseudo-holomorphic curves.
\end{abstract}
\maketitle
\date{}
\footnotetext[1]{Both authors are supported in part by NSF}

Recently, the quantum cohomology rings have generated 
a lot of researches. 
Many mathematicians have contributed to this rapidly progressing
field of mathematics. 
We will not described the history here,
but refer the interested reader to the orignal papers
and surveys (e.g.,
\cite{Beh}, \cite{Beh-Fan}, \cite{Beh-Man}, \cite{Fuk-Ono}, \cite{Giv-Kim}, 
\cite{Kon-Man1}, \cite{Kon-Man2},
\cite{Ler-Vaf-War}--\cite{Lu},
\cite{McD-Sal}--\cite{Vaf}, \cite{Wit1}, \cite{Wit2}
and the references therein).

The purpose of this paper is to give the
construction of  another deformation of the de Rham
cohomology ring.
The existence of a different deformation should not be a surprise,
since there is no reason to expect the deformation
of the cohomology to be unique.
A remarkable feature of our construction is that it
follows the traditional construction of the de Rham cohomology.
More precisely, we construct a quantum wedge product $\wh$
on exterior forms, and 
a quantum exterior differential $d_h$,
which satisfy the usual property of the calculus of 
differential forms.
This quantum calculus allows us  to 
``deformation quantize" many differential
geometric  objects, i.e. our quantum objects is a polynomial in
an indeterminate $h$, 
whose zeroth order terms are the classical objects. 
(In this sense, $h$ should be regarded as the Planck constant.)
For example, we will define quantum curvature of an ordinary
connection, and define quantum characteristic classes 
in the same fashion as the classical Chern-Weil theory.
 Our construction  has the following features which are not
shared by the quantum cohomology:
\begin{enumerate}
\item Quantum de Rham cohomology can be defined for 
Poisson manifolds,
not necessarily compact, or closed.
\item The proof of associative is of elementary nature.
\item It is routine to define quantum Dolbeault cohomology.
\item It is routine to define quantum characteristic class.
\item It is routine to define quantum equivariant 
de Rham cohomology.
\item The computations for homogeneous examples are elementary.
\end{enumerate}

Our construction is motivated by Moyal-Weyl multiplication
and Clifford multiplication.
For any finite dimensional vector space $V$ with a basis 
$\{ e_1, \cdots, e_m \}$,
let $\{ e^1, \cdots, e^m \}$ be the dual basis.
Assume that $w = w^{ij} e_i \otimes e_j \in V \otimes V$,
then $w$ defines a multiplication $\wedge_w$ on $\Lambda(V^*)$,
and a multiplication $*_w$ on $S(V^*)$,
such that $e^i \wedge_w e^j =   e^i \wedge e^j + w^{ij}$,
$e^i *_w e^j = e^i \odot  e^j + w^{ij}$.
If $w \in S^2(V)$, then $\wedge_w$ is the Clifford multiplication.
If $w \in \Lambda^2(V)$ is nondegenerate, $*_w$ is 
the Moyal-Weyl multiplication.
If $w \in \Lambda(V)$, then $\wedge_w$ is what we call a quantum 
exterior product (or a quantum Clifford multiplication).
It is elementary to show that this mutiplication is associative. 
We will use it to obtain a  quantum calculus on any Poisson manifold.
The main results we obtained in this paper have been announced in 
\cite{cao-zho}.

The layout of this paper is clear from the following

\tableofcontents

\noindent {\bf Acknowledgement} 
{\em The work in this paper was carried out during the second 
author's visit at Texas A\&M University.
He thanks the Department of Mathematics, especially the 
Geometry-Analysis-Topology group for hospitality and financial support.
Both authors want to thank Tony Giaquinto for a wonderful talk
on Poisson-Lie groups which jump started our initially unsuccessful approach.
We also thank him for providing Vaisman's book \cite{Vai},
which is a very good reference.}

\medskip
\begin{center}
{\bf \large Part I. Algebraic Theory} 
\end{center}

\section{Quantum Exterior Algebra}

\subsection{Deformation quantization} \label{sec:deform}

For more information deformation quantization of algebras, 
we refer to Donin \cite{Don} and the references therein.
Let $A$ be an algebra with unit over a field $\bk$ 
of characteristic zero.
A {\em deformation quantization} of $A$ 
is an algebra $A_h$ over $\bk[h]$ that is isomorphic to 
$A[h] = A \otimes_{\bk} \bk[h]$ as a $\bk[h]$-module,
such that $A_h/hA_h \cong A$.
A deformation quantization of  an algebra $A$ is uniquely 
determined by a $\bk$-linear map 
$f: A \otimes_{\bk} A \rightarrow A[h]$, 
$f(a, b) = a \cdot b + \sum_{j > 0} f_j(a, b) h^j$, 
for $a, b \in A$, 
where $a \cdot b$ stands for multiplication in $A$.
When $A$ is a ${\Bbb Z}$-graded algebra, 
a {\em graded deformation quantization} of $A$ is 
a deformation quantization $A_h \cong A \otimes_{\bk} \bk[h]$, 
which has the structure of a graded algebra when 
it is given the induced grading by setting $\deg (h) = 2$.

Assume now $A$ is a graded differential algebra (GDA), i.e,
there is a ${\Bbb Z}$-grading on $A$, and a $\bk$-linear 
$d: A \rightarrow A$ is a derivation of degree $1$ on $A$, such that 
$d^2 = 0$, and  $d(a\cdot b) = (da) \cdot b + (-1)^{|a|} a \cdot (db)$,
for homogeneous elements $a, b \in A$, 
and $|a|$ stands for the degree of $a$.
The graded algebra $H^*(A, d) = \Ker d /\Img d$ is called the 
cohomology of the GDA $(A, d)$. 
A deformation quantization of a GDA $(A, d)$ is a 
graded deformation quantization $A_h$, 
together with a $\bk[h]$-differential $d_h$ of degree $1$, 
such that when $A_h/hA_h$ is identified with $A$, 
the map on $A_h/hA_h$ induced by
$d_h$ is identified with $d$.
In this paper, we will be concerned with $d_h$  of the form
$d_h = d - h\delta$, where $\delta: A \rightarrow A$ is a derivation of
degree $-1$. 
We regard the complex $(A[h], d_h)$ as 
as the associated complex of the double complex 
$(C^{p, q}: = h^pA^{q-p}, d, -h\delta)$.
There are two spectral sequences associated with 
this double complex by standard theory (Bott-Tu \cite{Bot-Tu}), 
one of them has $E_1^{p, q} = h^pH^{q-p}(A, d)$.
If this spectral sequence degenerate at $E_1$, we then have
$$Q_hH^n(A) : = H^n(A_h, d_h) = 
\oplus_{p+q = n} h^pH^{q-p}(A, d),$$
for $p, q \geq 0$. 
It follows then $Q_hH^*(A)$ is a deformation quantization 
of $H^*(A, d)$ in this case.

\subsection{Laurent deformation}
Sometimes it is useful to use 
${\bk}[h, h^{-1}]$ instead of ${\bk}[h]$ as the coefficient ring. 
This will become apparent in our theory in $\S \ref{sec:dRCoh}$.
Let $A$ be an algebra with unit over a field of characteristic zero,
a (polynomial) Laurent deformation of $A$ is an algebra over $\bk[h, h^{-1}]$,
which is isomorphic to $A[h, h^{-1}] = A \otimes_{\bk} \bk[h, h^{-1}]$,
whose multiplication is determined by a $\bk$-linear map
$$f: A \otimes_{\bk} A \rightarrow A[h, h^{-1}]$$
of the following form: $f(a, b) = a \cdot b + \sum_{j\neq 0}f_j(a, b)h^j$,
where $a, b \in A$, and $a\cdot b$ stands for multiplication in $A$.
We use the following simple construction.
Since  a deformation quantization $(A_h, *_h)$ of an algebra $(A, \cdot)$
is  determined by a $\bk$-linear map
$f: A \otimes_{\bk} A \rightarrow A[h]$ of the form $f(a, b) = a \cdot b
+ \sum_{j > 0} f_j(a, b)h^j$, 
it gives rise to a unique Laurent deformation of $A$.
We will consider Laurent deformation of a GDA similar to the polynomial 
deformation case discussed in $\S \ref{sec:deform}$.

\subsection{Moyal-Weyl quantization}

Reformulation and generalization of
Moyal-Weyl quantization for  polynomial algebra on
a symplectic vector space is presented in this section.
It serves as a motivation for our construction below.

Let $V$ be a $\bk$-vector space. 
The symmetric tensor algebra $S(V^*)$ of 
$V^*$ can be regarded as the algebra of $\bk$-polynomials on $V$.
Let $w \in \Lambda^2(V)$. 
When a basis $\{e_1, \cdots, e_n\}$ is chosen,
let $(x^1, \cdots, x^n)$ be the coordinates in this basis,
$w = w^{ij}e_i \wedge e_j$.
Then the Moyal-Weyl product of two polynomials $u, v$ is given by
\begin{eqnarray*} 
u*_h v = \sum_{n \geq 0} \frac{h^n}{n!} 
	w^{i_1j_1}\cdots w^{i_nj_n}
	\frac{\partial^n u}{\partial x^{i_1}\cdots \partial x^{i_n}} \cdot
	\frac{\partial^n v}{\partial x^{j_1}\cdots \partial x^{j_n}}.
\end{eqnarray*}
It can be formulated without using coordinates.
Denote by $T(V^*)$ the tensor algebra of $V^*$, 
for any  $\phi \in \otimes^2V$,
define $L_{\phi}: T(V^*) \otimes T(V^*) \rightarrow T(V^*)$ as follows:
for $\alpha, \beta \in T(V^*)$, $L_{\phi}(\alpha \otimes \beta)$
is obtained from $\alpha \otimes \phi \otimes \beta$
by contracting the first factor in $\phi$ with the last
factor in $\alpha$, and contraction the second factor in $\phi$ with
the first factor in $\beta$.
In coordinates, this is denoted by
$$L_{\phi}(\alpha \otimes \beta) 
	= \phi^{ij} (\alpha \vdash e_i) \otimes 
	(e_j \dashv \beta),$$
where $\vdash$ is the contraction from the right, 
$\dashv$ is the contraction from the left,
i.e.,
\begin{eqnarray*}
(\alpha \vdash v)(v_1, \cdots, v_{k-1}) 
	= \alpha(v_1, \cdots, v_{k-1}, v), \\
(v \dashv \alpha)(v_1, \cdots, v_{k-1}) 
	= \alpha(v, v_1, \cdots, v_{k-1}),
\end{eqnarray*}
for $v, v_1, \cdots, v_{k-1} \in V$, $\alpha \in T^{k}(V^*)$.
Denote by $m_s: S(V^*) \otimes S(V^*) \rightarrow S(V^*)$ the
symmetric product.	
It is then obvious that Moyal-Weyl product is given by
$$\alpha *_h \beta = m_s( \exp (hL_w) (\alpha \otimes \beta)).$$
We give $S_h(V^*) = S(V^*)[h]$ the following $\Bbb Z$-grading: 
elements in $S^p(V^*)$ has degree $p$, and $h$ has degree $2$.
If we denote by $S_h^m(V^*)$ the space of homogeneous elements
of degree $m$, then we have 
$$S_h^{[m]}(V^*) *_h S_h^{[n]}(V^*) \subset S_h^{[m+n]}(V^*).$$
Since $w$ is anti-symmetric, $u*_h v = v *_h u$ does not hold
in general.
It can be checked that $(u*_hv)*_hw = u*_h(v*_hw)$.
See e.g Bayen-Flato-Fronsdal-Lichnerowicz-Sternheimer 
\cite{Bay-Fla-Fro-Lic-Ste}
or the next section.
Therefore, $(S_h(V^*), *_h)$ is an algebra.
It is clear that it is a deformation quantization of 
the polynomial algebra $S(V^*)$.

\subsection{Quantum exterior algebra}

Now let $\Lambda(V^*)$ denote the exterior algebra, and 
$m: \Lambda(V^*) \rightarrow \Lambda(V^*)$ the exterior product.
Given any $w = w^{ij}e_i  \wedge e_j \in   \Lambda^2(V)$,
we define the {\em quantum exterior product} by
$\wh: \Lambda(V^*) \otimes \Lambda(V^*) \rightarrow 
	\Lambda(V^*)[h]$ by
\begin{eqnarray*}
&   & \alpha \wedge_{h, w} \beta 
	= m(\exp (hL_w) (\alpha \otimes \beta))
	= \sum_{n \geq 0} \frac{h^n}{n!} m(L_{w}^n(\alpha \otimes \beta)) \\
& = & \sum_{n \geq 0} 
	\frac{h^n}{n!} w^{i_1j_1} \cdots w^{i_nj_n}
	(\alpha \dashv e_{i_1} \dashv \cdots \dashv e_{i_n})
	\wedge (e_{j_n} \vdash \cdots \vdash e_{j_1} \vdash \beta), 
\end{eqnarray*}	
for $\alpha, \beta \in \Lambda(V^*)$.
(This is evidenly independent of the choice of the basis.) 
Notice that this is just Moyal-Weyl multiplication for exterior algebra.
It is defined this way to keep track of the signs.
When there is no confusion about $w$, we will simply write
$\alpha \wh \beta$ for $\alpha \wedge_{h, w} \beta$. 
The map $\Lambda(V^*) \otimes \Lambda(V^*) \rightarrow \Lambda(V^*)$
will be denoted by $\ww$. We have
\begin{eqnarray*}
&   & \alpha \wh \beta	\\
& = & \sum_{n \geq 0} \frac{h^n}{n!}(-1)^{n(|\alpha| -1)}
	w^{i_1j_1} \cdots w^{i_nj_n}
	(e_{i_1} \vdash \cdots \vdash e_{i_n} \vdash  \alpha)
	\wedge (e_{j_n} \vdash \cdots \vdash e_{j_1} \vdash \beta) \\
& = & \sum_{n \geq 0} \frac{h^n}{n!}(-1)^{n(|\alpha| -1) + n(n-1)/2}
	w^{i_1j_1} \cdots w^{i_nj_n} \\
&   & \cdot	(e_{i_1} \vdash \cdots \vdash e_{i_n}\vdash  \alpha)
	\wedge (e_{j_1} \vdash \cdots \vdash e_{j_n} \vdash \beta).
\end{eqnarray*}
For simplicity of the notations, we will write
\begin{eqnarray*}
\alpha \wh \beta = 
	\sum_{n \geq 0} \frac{h^n}{n!}(-1)^{n(|\alpha| -1) + n(n-1)/2}
	\sum_{\substack{|I|=n \\ |J| = n}}
	w^{IJ} (e_I\vdash \alpha)\wedge (e_J \vdash \beta).
\end{eqnarray*}
We extend $\wh$ as $\bk[h]$-module map to 
$\Lambda_h(V^*) \otimes_{\bk[h]} \Lambda_h(V^*)$.
We give the same $\Bbb Z$-grading on 
$\Lambda_h(V^*) = \Lambda(V^*)[h]$
 as for $S_h(V^*)$, then it is clear that
$$\Lambda_h^{[m]}(V^*) \wedge_h \Lambda_h^{[n]}(V^*) \subset
	\Lambda_h^{[m+n]}(V^*).
$$

\begin{remark} \label{Clifford}
To make contact with more familiar objects, we consider the case of
$|\alpha| = 1$, we take $\alpha$ to be an element $v \in V^*$.
The bivector $w \in \Lambda^2(V^*)$ defines a homomorphism
$V^* \rightarrow V$ by $v \mapsto v_{w} = w(v, \cdot)$.
Then $v \wh \beta = v \wedge \beta + h {v_{w}} \vdash \beta$.
This is the analogue of 
the Clifford multiplication, which is defined by setting $h = 1$, 
and using an element $q \in S^2(V)$ instead of $w$. 
So it might be more instructive to call $\wh$ the 
quantum skew Clifford multiplication.
Of course, when $h = 0$, both quantum Clifford multiplication 
and quantum skew Clifford multiplication gives the exterior product.
This reveals that both of them are deformation quantizations of 
the exterior product, the difference being one preserves 
the super commutativity, the other destroys it.
\end{remark}

\begin{theorem} \label{thm:ass}
The quantum exterior product satisfies the following properties
\begin{eqnarray}
\text{Supercommutativity} 
& & \alpha \wh \beta = (-1)^{|\alpha||\beta|} \beta \wh \alpha, 
	\label{com} \\
\text{Associativity} 
& & (\alpha \wh \beta) \wh \gamma = \alpha \wh (\beta \wh \gamma), 
	\label{ass}
\end{eqnarray}
for  all $\alpha, \beta, \gamma \in \Lambda_h(V^*)$.
Therefore, $(\Lambda_h(V^*), \lambda_h)$ is a deformation quantization
of the exterior algebra $(\Lambda(V^*), \wedge)$.
\end{theorem}

\begin{remark}
Similar results hold for $\Lambda_{h, h^{-1}}(V^*)$, and for $\ww$
on $\Lambda(V^*)$.
\end{remark}

\begin{proof}[Proof of supercommutativity]
\begin{eqnarray*}
&    & \beta \wh \alpha = 
	\sum_{n \geq 0} \frac{h^n}{n!}(-1)^{n(|\beta| -1) + n(n-1)/2}
	\sum_{\substack{|I|=n \\ |J| = n}}
	w^{IJ} (e_I\vdash \beta)\wedge (e_I \vdash \alpha) \\
&  = & \sum_{n \geq 0} \frac{h^n}{n!}(-1)^{n(|\beta| -1) + n(n-1)/2}
	(-1)^n \sum_{\substack{|I|=n \\ |J| = n}}w^{JI} 
	(-1)^{(|\alpha|-n)(|\beta|-n)}
	(e_J\vdash \alpha)\wedge (e_I \vdash \beta) \\
& = & (-1)^{|\alpha||\beta|} \sum_{n \geq 0} \frac{h^n}{n!}(-1)^{n(|\alpha| -1) + n(n-1)/2}
	\sum_{\substack{|I|=n \\ |J| = n}}
	w^{IJ} (e_I\vdash \alpha)\wedge (e_J \vdash \beta) \\
& = & (-1)^{|\alpha||\beta|}\beta \wh \alpha.
\end{eqnarray*}
\end{proof}

We will prove the associativity $(\ref{ass})$ by induction.
We say that $A(a, b,c )$ holds, if $(\ref{ass})$ holds 
for all $\alpha, \beta, \gamma \in \Lambda(V^*)$ with
$|\alpha| = a$, $|\beta| = b$, $|\gamma| = c$. 
We say that $A(\leq a, b, c)$ holds, if  $(\ref{ass})$ holds 
for all $\alpha, \beta, \gamma \in \Lambda(V^*)$ with
$|\alpha| \leq a$, $|\beta| = b$, $|\gamma| = c$.
Our strategy is as follows.
We first prove $A(1, b, c)$ for arbitrary $b$ and $c$, then 
prove the general case by induction on $a$.

\begin{proof}[Proof of $A(1, b, c)$]
By linearity, we can assume that  $\alpha = e^i$ 
for some given basis $\{e_1, \cdots, e_m\}$ of $V$. 
\begin{align}
&  (e^i \wh \beta) \wh \gamma \nonumber \\
 =& (e^i \wedge \beta + h w^{ik}e_k \vdash \beta) 
 	\wh \gamma \nonumber\\
 = & \sum_{n \geq 0} \frac{h^n}{n!}(-1)^{n|\beta| + n(n-1)/2}
	\sum_{\substack{|I|=n \\ |J| = n}}
	w^{IJ}(e_I \vdash (e^i \wedge \beta))
	\wedge (e_J \vdash \gamma)  \tag{\bf A}  \\
& +   \sum_{n \geq 0} \frac{h^{n+1}}{n!}(-1)^{n|\beta| + n(n-1)/2}
	\sum_{\substack{|I|=n \\ |J| = n}}
	w^{ik}w^{IJ}
	(e_I \vdash e_k \vdash \beta)
	\wedge (e_J \vdash \gamma) \tag{\bf B}
\end{align}
Now  for $(\bf B)$, we use 
$e_I \vdash e_k \vdash \beta = (-1)^n e_k \vdash e_I \vdash \beta$.
For $(A)$, we use 
$$
e_I \vdash (e^i \wedge \beta) =   
(-1)^n e^i \wedge (e_I \vdash \beta) 	
+ \sum_{l =1}^n (-1)^{n-l}\delta^{i}_{i_l} 
	(e_{i_1} \vdash \cdots \vdash \hat{e}_{i_l} \vdash \cdots \vdash 
		e_{i_n} \vdash \beta).
$$
This still holds when $n = 0$, 
if we take the second term on the right to be zero.
So $({\bf A})$ is equal to 
\begin{align}
 &  \sum_{n \geq 0} \frac{h^n}{n!}(-1)^{n|\beta| + n(n-1)/2}
 	\sum_{\substack{|I|=n \\ |J| = n}}
 	w^{IJ} (-1)^n 
 	 e^i \wedge (e_I \vdash \beta)
	\wedge (e_J \vdash \gamma) \tag{$\bf A_1$} \\
& + \sum_{n \geq 1} \frac{h^n}{n!}(-1)^{n|\beta| + n(n-1)/2}
	\sum_{\substack{|I|=n \\ |J| = n}} 
	w^{IJ} \tag{$\bf A_2$} \\
&	\cdot \sum_{1\leq l \leq n} (-1)^{n-l}\delta^{i}_{i_l} 
	(e_{i_1} \vdash \cdots \vdash \hat{e}_{i_l} \vdash \cdots \vdash 
		e_{i_n} \vdash \beta)
	\wedge (e_J \vdash \gamma) \nonumber
\end{align}
For $(\bf A_2)$, we use the following  renaming of the indices:
$I'_l = i_1 \cdots \hat{i}_l \cdots i_n$, 
$J'_l= j_1 \cdots \hat{j}_l \cdots j_n$, $j_l = k$. 
Then we have
\begin{eqnarray*}
 & & \sum_{\substack{|I|=n \\ |J| = n}} w^{IJ} 
	 \sum_{1\leq l \leq n} (-1)^{n-l}\delta^{i}_{i_l} 
	(e_{i_1} \vdash \cdots \vdash \hat{e}_{i_l} \vdash \cdots \vdash 
		e_{i_1} \vdash \beta) \wedge (e_J \vdash \gamma) \\
& = & \sum_{l=1}^n \sum_{\substack{|I|=n \\ |J| =n}} 
(-1)^{n-l}w^{ik} w^{I'_lJ'_l} 
(e_{I'_l} \vdash \beta) 
	\wedge (-1)^{l-1} (e_k \vdash e_{J_l'} \wedge \gamma) \\
& = & n \sum_{\substack{|I|=n \\ |J| =n-1}} 
	(-1)^{n-1} w^{ik} w^{IJ} (e_I \vdash \beta) 
	\wedge (e_k \vdash e_J \wedge \gamma).
\end{eqnarray*}
So $(\bf A_2)$ is equal to 
\begin{eqnarray*}
&   & \sum_{n \geq 1} \frac{h^n}{n!}(-1)^{n|\beta| + n(n-1)/2}
	n (-1)^{n-1} \sum_{\substack{|I|=n \\ |J| =n-1}}
	w^{ik} w^{IJ} (e_I \vdash \beta) 
	\wedge (e_k \vdash e_J \vdash \gamma) \\
& = & \sum_{n \geq 1} \frac{h^n}{(n-1)!}
	(-1)^{n|\beta| + (n-2)(n-1)/2}
 	\sum_{\substack{|I|=n \\ |J| =n-1}}
 	w^{ik} w^{IJ} (e_I \vdash \beta) 
	\wedge (e_k \vdash e_J \vdash \gamma) \\
& = & \sum_{n \geq 0} \frac{h^{n+1}}{n!}
	(-1)^{(n+1)|\beta| + n(n-1)/2}
	\sum_{\substack{|I|=n \\ |J| =n}}
 	w^{ik} w^{IJ} (e_I \vdash \beta) 
	\wedge (e_k \vdash e_J \vdash \gamma) 
\end{eqnarray*}
To summarize, we have
\begin{align}
&    (e^i \wh \beta) \wh \gamma  
	\nonumber \\
& = \sum_{n \geq 0} \frac{h^n}{n!}(-1)^{n|\beta| + n(n+1)/2}
	\sum_{\substack{|I|=n \\ |J| =n}} w^{IJ} 
 	 e^i \wedge (e_I \vdash\beta)
	\wedge (e_J \vdash \gamma) \tag{$\bf A_1$} \\
& + \sum_{n \geq 0} \frac{h^{n+1}}{n!}(-1)^{(n+1)|\beta| + n(n-1)/2}
 	\sum_{\substack{|I|=n \\ |J| =n}}
 	w^{ik} w^{IJ} (e_I \vdash \beta) 
	\wedge (e_k \vdash e_J \vdash \gamma) \tag{$\bf A_2$} \\
& + \sum_{n \geq 0} \frac{h^{n+1}}{n!}(-1)^{n|\beta| + n(n+1)/2}
	\sum_{\substack{|I|=n \\ |J| =n}}
	w^{ik} w^{IJ}
	(e_k \vdash e_I \vdash \beta)
	\wedge (e_J \vdash \gamma) \tag{\bf B}
\end{align}
Similarly,
\begin{align}
&  e^i \wh (\beta \wh \gamma) \nonumber \\
 =& e^i \wh \sum_{n \geq 0} \frac{h^n}{n!}(-1)^{n(|\beta|-1) + n(n-1)/2}
	\sum_{\substack{|I|=n \\ |J| =n}}w^{IJ}
	(e_I \vdash \beta)
	\wedge (e_J \vdash \gamma)  \nonumber \\
 = &  e^i  \wedge \sum_{n \geq 0} \frac{h^n}{n!}(-1)^{n(|\beta|-1) + n(n-1)/2}
	\sum_{\substack{|I|=n \\ |J| =n}}w^{IJ}
	(e_I \vdash  \beta)
	\wedge (e_J \vdash \gamma)  \tag{\bf C} \\
& +  h \sum_{n \geq 0} \frac{h^n}{n!}(-1)^{n(|\beta|-1) + n(n-1)/2}
	\sum_{\substack{|I|=n \\ |J| =n}}w^{ik} w^{IJ}
	e_k\vdash [(e_I \vdash  \beta)
	\wedge (e_J \vdash \gamma)]  \tag{\bf D} \\
 = &  	e^i  \wedge \sum_{n \geq 0} \frac{h^n}{n!}(-1)^{n|\beta| + n(n+1)/2}
	\sum_{\substack{|I|=n \\ |J| =n}}w^{IJ}
	(e_I \vdash  \beta)
	\wedge (e_J \vdash \gamma)  \tag{\bf C} \\
& +    \sum_{n \geq 0} \frac{h^{n+1}}{n!}(-1)^{n|\beta| + n(n+1)/2}
	\sum_{\substack{|I|=n \\ |J| =n}}w^{ik} w^{IJ}
	(e_k\vdash e_I \vdash  \beta)
	\wedge (e_J \vdash \gamma)]  \tag{$\bf D_1$} \\
& +    \sum_{n \geq 0} \frac{h^{n+1}}{n!}(-1)^{n|\beta| + n(n+1)/2} 
	\tag{$\bf D_2$}\\
&	\cdot 
	\sum_{\substack{|I|=n \\ |J| =n}}w^{ik} w^{IJ}
	(-1)^{|\beta|-n} (e_I \vdash  \beta)
	\wedge (e_k \vdash e_J \vdash \gamma)]  
	\nonumber
\end{align}
It is clear that $({\bf A_1}) = ({\bf C})$, $({\bf A_2}) = ({\bf D_2})$,
and $({\bf B}) = ({\bf D_1})$. 
This completes the proof of $A(1, b, c)$.
\end{proof}

\begin{proof}[Proof of $A(a, b, c)$]
Assume that $A(\leq a, b,c)$ is proven, 
we now show how to deduce $A(1+ a, b, c)$.
Without loss of generality, assume that $\alpha = v \wedge \eta$,
for some $v$ and $\eta$ with $|v| = 1$, $|\eta| = a$.
By definition, $v \wh \eta = v \wedge \eta - h f(v, \eta)$
for some element $f(v, \eta)$ with degree $\leq a - 2$.
Therefore,
\begin{align*}
& \alpha \wh (\beta \wh \gamma)
 = (v \wedge \eta) \wh (\beta \wh \gamma) & \\
 = &(v \wh \eta + hf(v, \eta)) \wh (\beta \wh \gamma) & \\
 = & (v \wh \eta) \wh (\beta \wh \gamma) + hf(v, \eta) \wh (\beta \wh \gamma) & \\
 = &v \wh (\eta \wh (\beta \wh \gamma)) + h(f(v, \eta) \wh \beta) \wh \gamma & 
	\text{by $A(1, a, b+c)$ and  $A(\leq a, b, c)$} \\
 = &v \wh ((\eta \wh \beta) \wh \gamma) + h(f(v, \eta) \wh \beta) \wh \gamma & 
	\text{by $A(a, b, c)$} \\
 = & (v \wh (\eta \wh \beta)) \wh \gamma + h(f(v, \eta) \wh \beta) \wh \gamma & 
	\text{by $A(1, a+b, c)$} \\
 = & ((v \wh \eta) \wh \beta) \wh \gamma + h(f(v, \eta) \wh \beta) \wh \gamma & 
	\text{by $A(1, a, b)$} \\
 = & ((v \wh \eta + hf(v, \eta)) \wh \beta) \wh \gamma  &\\
 = & (v \wedge \eta) \wh \beta) \wh \gamma = (\alpha \wh \beta) \wh \gamma. &
\end{align*}
Therefore, associativity 
$(\ref{ass})$ holds for all 
$\alpha, \beta, \gamma  \in \Lambda_h(V^*)$. 
\end{proof}

Two remarks are in order. 
First, in the above proof of associativity,
we have never used the anti-symmetric property of $w$.
Therefore, if we define $\wh$ using any bi-vector $\phi$, 
$\wh$ will still be associative. 
It may not be supercommutative anymore.
Second, the associativity of Moyal-Weyl quantization 
can be proved in the same fashion. 
Without those $\pm$ signs, it is much simpler.
Again we do not use the anti-symmetric property of $w$, so
if one defines generalized Moyal-Weyl quantization $*_h$ using
any bi-vector $\phi$, $*_h$ is associative.
In particular, if $\phi$ is symmetric, $*_h$ defined by $\phi$ is
commutative.

\subsection{Complexified quantum exterior algebra}

In this section, we will consider real vector space $V$ with an
almost complex structure $J$, i.e., $J : V \rightarrow V$
is a linear transformation such that $J^2 = - Id$.
There is an induced linear transformation 
$\Lambda^2J: \Lambda^2(V) \rightarrow \Lambda^2(V)$.
For any bi-vector $w \in \Lambda^2(V)$,
we say  $J$ preserves $w$ if 
$\Lambda^2J (w) = w$.
Given any bi-vector $w$ which is preserved by $J$, 
we can define the quantum exterior product on 
$\Lambda_h(V^*)$ as in the last section.
Now if we tensor everything by $\Bbb C$,
we get a complex algebra ${\Bbb C}\Lambda_h(V^*)$, which is
a deformation quantization of ${\Bbb C}\Lambda(V^*) :=
\Lambda(V^*) \otimes_{\Bbb R} {\Bbb C} 
= \Lambda_{\Bbb C}(V^* \otimes_{\Bbb R} {\Bbb C})$.
As a common practice in complex geometry, 
we can exploit a natural decomposition as follows.
There are two complex vector spaces $V^{1, 0}$ and $V^{0, 1}$
with underlying real vector space  $V$ by:
for $V^{1, 0}$,  the multiplication by $\sqrt{-1}$ is given by $J$;
for $V^{0, 1}$, by $-J$.
There is a natural identification of complex vector
spaces  ${\Bbb C}V \cong V^{1, 0} \oplus V^{0, 1}$ given 
by $v = \frac{1}{2}(v - \sqrt{-1}Jv) + \frac{1}{2}(v - \sqrt{-1}Jv)$
for any $v \in V$, and extend it complex linearly to ${\Bbb C}V$.
As a consequence, there are decompositions
\begin{eqnarray*}
{\Bbb C}\Lambda(V) & = & \oplus_{p, q} \Lambda^{p, q}(V), \\
{\Bbb C}\Lambda(V^*) & = & \oplus_{p, q} \Lambda^{p, q}(V^*),
\end{eqnarray*}
where $\Lambda^{p, q}(V) \cong \Lambda_{\Bbb C}^p(V^{1, 0})
\otimes_{\Bbb C} \Lambda^q_{\Bbb C}(V^{0, 1})$,
and $\Lambda^{p, q}(V^*) \cong \Lambda_{\Bbb C}^p((V^{1, 0})^*)
\otimes_{\Bbb C} \Lambda^q_{\Bbb C}((V^{0, 1})^*)$.
We give ${\Bbb C}\Lambda_h(V^*)$ the following 
${\Bbb Z} \times {\Bbb Z}$-bigrading: 
elements in $\Lambda^{p, q}(V^*)$
has bi-degree $(p, q)$,  $h$ has bi-degree $(1, 1)$.
Since $w$ is preserved by $J$, it belongs to $\Lambda^{1, 1}(V)$
after complexification.
Denote by $\Lambda_h^{[p, q]}(V^*)$ the space of homogeneous
elements of bi-degree $(p, q)$.
When we compute $L$ in ${\Bbb C}\Lambda(V^*)$,
we can use a complex basis of ${\Bbb C}V$ of the form
$\{ f_1, \cdots, f_n, f_{\bar 1}, \cdots, f_{\bar n}\}$,
where $\{ f_1, \cdots, f_n \}$ is a complex basis of
$V^{1, 0}$, and $\{ f_{\bar 1}, \cdots, f_{\bar n}\}$ 
is the complex conjugate basis of $V^{0, 1}$. 
It is then clear from the definition that
\begin{eqnarray*}
\alpha \wh \beta & = & \sum_{p, q \geq 0} 
	\frac{h^{p+q}}{p!q!}
	\sum_{\substack{|A| = |B| = p \\ |C| = |D| = q}}
	w^{A\bar{B}}w^{\bar{C}D}
	(\alpha \dashv f_A \dashv f_{\bar{C}}) 
	\wedge
	(f_{r(D)} \vdash f_{r(\bar{B})} \vdash \beta),
\end{eqnarray*}
where if $A = (a_1, \cdots, a_p)$, then $r(A) = (a_p, \cdots, a_1)$,
the reverse of $A$. 
It then follows that
$$\Lambda_h^{[p, q]}(V^*) \wh \Lambda_h^{[r, s]}(V^*)
	\subset \Lambda_h^{[p+r, q+t]}(V^*).$$
	
Now let  $\omega$ be a symplectic 
form on $V$, which is compatible with an almost complex structure
on $V$, i.e. rank of $\omega$ is $2n = \dim (V)$, 
$w(J\cdot, J\cdot) = \omega(\cdot, \cdot)$,
and $g(\cdot, \cdot) := \omega(\cdot, J\cdot)$ is 
a positive definite element of $S^2(V^*)$.
For $X, Y \in  {\Bbb C}V$, set
$$H(X, Y) = \frac{1}{\sqrt{-1}}w(X, \overline{Y}).$$
Then $H$ is a Hermitian metric on ${\Bbb C}V$, such that
$V^{1, 0} \perp V^{0, 1}$.
It induces a Hermitian metric on ${\Bbb C}\Lambda(V^*)$,
which we will give explicitly in coordinates below.
It is possible to find an orthonormal basis of $V$ for $g$ of the 
form $\{e_1, Je_1, \cdots, e_n, Je_n\}$.
Set 
\begin{eqnarray*}
f_a  =  \frac{1}{2}(e_a - \sqrt{-1}Je_a), &
f_{\bar{a}}  =  \frac{1}{2}(e_a + \sqrt{-1}Je_a),
\end{eqnarray*}
for $a = 1, \cdots, n$.
It can be easily checked that
\begin{eqnarray*}
\omega_{ab} := \omega(f_a, f_b) = 0, \\
\omega_{\bar{a}\bar{b}} := \omega(f_{\bar{a}}, \omega_{\bar{b}}) = 0, \\
\omega_{a \bar{b}} = - \omega_{\bar{b}a} 
	= \frac{\sqrt{-1}}{2} \delta_{ab},
\end{eqnarray*}
for $a, b = 1, \cdots, n$. 
Then $\{\sqrt{2}f_a \}$ is an orthonormal basis of $V^{1, 0}$, and
$\{ \sqrt{2}f_{\bar a}\}$ is an orthonormal basis of $V^{0, 1}$.
Let $\{\frac{1}{\sqrt{2}}f^a \}$ and 
$\{ \frac{1}{\sqrt{2}}f^{\bar a}\}$ be the dual
basis for $(V^{1, 0})^*$ and $(V^{0, 1})^*$ respectively.
Then in the canonically induced Hermitian metric 
on $\Lambda^{p, q}(V^*)$, 
$\{ 2^{-(p+q)/2}f^{a_1} \wedge \cdots \wedge f^{a_p} \wedge
f^{\bar{b}_1} \wedge \cdots \wedge f^{\bar{b}_q}, 
a_1 < \cdots < a_p, b_1 < \cdots < b_q \}$
is an orthonormal basis.
The symplectic form $\omega \in \Lambda^2(V^*)$ determines 
a unique bi-vector 
$w^{\sharp} \in \Lambda^2(V)$ 
in a way similar to raising the index in Riemannian geometry.
Let $\omega^{ij} = \omega^{\sharp}(e^i, e^j)$, 
$\omega_{ij} = \omega(e_i, e_j)$, then 
the matrix $(\omega^{ij})$ is the 
inverse of the matrix $(\omega_{ij})$.
After complexification, let $a, b$ denotes the complex indices,
then we have
\begin{eqnarray*}
\omega^{ab} = \omega^{\bar{a}\bar{b}} = 0, & 
\omega^{a\bar{b}} = - w^{\bar{b}a} 
	= -\frac{2}{\sqrt{-1}}\delta^{ab}.
\end{eqnarray*}
For any $\alpha \in \Lambda^{p, q}(V^*)$,
we can write 
$$\alpha = \sum_{\substack{|A| =p \\ |C| = q}} \frac{1}{p!q!}
\alpha_{A\bar{C}} f^A \wedge f^{\bar{C}},$$
where $\alpha_{A\bar{C}}$ is anti-symmetric in the complex
indices $A$ and $C$. 
Then we have
$$\alpha_{A\bar{C}} = f_{r(\bar{C})} \vdash f_{r(A)} \vdash \alpha.$$
Therefore, for $\alpha, \beta \in \Lambda^{p, q}(V^*)$,
\begin{eqnarray*}
&   &  H(\alpha, \beta) \\
& = &   \sum_{\substack{|A| = p \\|C| = q}} 
        \frac{2^{p+q}}{p!q!} 
        (f_{r(\bar{C})} \vdash f_{r(A)} \vdash \alpha)
	\cdot
	(f_{r(C)} \vdash f_{r(\bar{A})} \vdash \bar{\beta}) \\
& = & \sum_{\substack{|A| = |B| = p
	\\|C| = |D| = q}} \frac{2^{p+q}}{p!q!}
	\delta^{AB} \delta^{CD} 
	(f_{r(\bar{C})} \vdash f_{r(A)} \vdash \alpha)
	\cdot
	(f_{r(D)} \vdash f_{r(\bar{B})} \vdash \bar{\beta}) \\
& = & \sum_{\substack{|A| = |B| = p
	\\|C| = |D| = q}} \frac{(\sqrt{-1})^{p-q}}{p!q!}
	(-1)^p
	\omega^{A\bar{B}}\omega^{\bar{C}D}
	(f_{r(\bar{C})}\vdash f_{r(A)} \vdash   \alpha) 
	\cdot (f_{r(D)} \vdash f_{r(\bar{B})}\vdash \bar{\beta}) \\
& = & \frac{(\sqrt{-1})^{p-q}}{p!q!}
	(-1)^{p + (p+q)(p+q-1)/2} \\
&   &  \cdot	\sum_{\substack{|A| = |B| = p
	\\|C| = |D| = q}} 
	\omega^{A\bar{B}}\omega^{\bar{C}D}
	(\alpha \dashv f_A \dashv f_{\bar{C}}) 
	\cdot (f_{r(D)} \vdash f_{r(\bar{B})}\vdash \bar{\beta}) \\
& = & \frac{(\sqrt{-1})^{p-q}}{2^{p+q}}
	(-1)^{p + (p+q)(p+q-1)/2} (\alpha \ww \bar{\beta})_0,
\end{eqnarray*}
where $\alpha \ww \bar{\beta}$ is obtained 
from $\alpha \wh \bar{\beta}$ by
 setting $h = 1$. And the subscript $0$ means taking the degree
 $0$ zero part.
Since for general $\alpha \in \Lambda^{p, q}(V^*)$,
$\beta \in \Lambda^{r, s}(V^*)$, $(\alpha \ww \bar{\beta})_0$ can be
nonzero only if $p = r$ and $q =t$,
we see that for any $\alpha \in \Lambda^{p, q}(V^*)$,
$\beta \in \Lambda^{r, s}(V^*)$, we have
\begin{equation} \label{eqn:metric}
\begin{split} 
H(\alpha, \beta) & =  (\sqrt{-1})^{p-q}
	(-1)^{p + (p+q)(p+q-1)/2} (\alpha \ww \bar{\beta})_0 \\
 & =  (\sqrt{-1})^{r-s}
	(-1)^{r + (r+s)(r+s-1)/2} (\alpha \ww \bar{\beta})_0.
\end{split}
\end{equation}

\begin{lemma}
For any three elements $\alpha \in \Lambda^{p, q}(V^*)$, 
$\beta \in \Lambda^{s, t}(V^*)$
and $\gamma \in \Lambda^{u, v}(V^*)$,
we have
$$ H(\alpha \ww \beta, \gamma)  
= H(\alpha, \beta \ww \gamma).$$
\end{lemma}

\begin{proof} 
>From $(\ref{eqn:metric})$, we have
\begin{eqnarray*}
&   & H(\alpha \ww \beta, \gamma) \\
& = &	(\sqrt{-1})^{u-v}(-1)^{u + (u+v)(u+v-1)/2}
	((\alpha \ww \beta) \ww \bar{\gamma})_0 \\
& = &	(\sqrt{-1})^{u-v}(-1)^{u + (u+v)(u+v-1)/2}
	(\alpha \ww \beta \ww \bar{\gamma})_0.
\end{eqnarray*}
It is nonzero only if 
$p+s - u = q + t - v \geq 0$.
Similarly 
\begin{eqnarray*}
&   & H(\alpha, \bar{\beta} \ww \gamma) \\
& = &	(\sqrt{-1})^{p-q}(-1)^{p + (p+q)(p+q-1)/2}
	(\alpha \ww \overline{\bar{\beta} \ww \gamma})_0 \\
& = &	(\sqrt{-1})^{p-q}(-1)^{p + (p+q)(p+q-1)/2}
	(\alpha \ww \beta \ww \bar{\gamma})_0.
\end{eqnarray*}
It is nonzero only if
$u + s - p = v + t - q \geq 0$. 
By associativity of $\wh$, $ H(\alpha \ww \beta, \gamma)$
and $H(\alpha, \beta \ww \gamma)$ differ only by a constant factor.
So we need only to determine this factor when 
both are nonzero.
In this case, we must have $s = t$.
hence $u - p = v - q$.
The factor then is
\begin{eqnarray*}
&  & (\sqrt{-1})^{(u - v) - (p -q)} 
	(-1)^{(u - p) + (u+v)(u+v-1)/2 - (p+q)(p+q-1)/2} \\
& = & (-1)^{(u - p) + (u+v+p+q)(u+v-p -q)/2 -(u+v - p - q)/2} \\
& = & (-1)^{(u - p)+ (2(u-p) + 2(p+q))(u-p) - (u - p)} = 1.
\end{eqnarray*}

\end{proof}

\subsection{Multiparameter deformation}

Let $\vec{w} = (w_1, \cdots, w_m)$, $w_j \in \Lambda^2(V)$, 
$j =1, \cdots, m$,  
$\vec{h} = (h_1, \cdots, h_m)$, 
$\Lambda_{\vec{h}}(V^*)=\Lambda(V^*)[\vec{h}]
= \Lambda(V^*)[h_1, \cdots, h_m]$.
Define 
$\wedge_{\vec{h}}: \Lambda(V^*) \otimes \Lambda(V^*)
\rightarrow  \Lambda(V^*)[\vec{h}]$ by
\begin{eqnarray*}
& &\alpha \wedge_{\vec{h}, \vec{\omega}} \beta 
  = m( \exp (h_1L_{w_1} + \cdots + h_mL_{w_m}) (\alpha \otimes \beta))\\
& = & \sum_{n_1, \cdots, n_m} \frac{h_1^{n_1}\cdots h_m^{n_m}}
	{n_1! \cdots n_m!}
	\sum_{\substack{|I_j| = n_j \\ |J_j| = n_j}}
	w_1^{I_1J_1}\cdots w_m^{I_mJ_m} \\
&  & 	(\alpha \dashv  e_{I_1} \dashv \cdots \dashv e_{I_m}) \wedge
	(e_{r(J_m)} \vdash \cdots \vdash e_{r(J_1)} \vdash \beta),
\end{eqnarray*}

\begin{theorem}
For any $\vec{w}$ and $\vec{h}$ as above,
$\wedge{\vec{h}, \vec{w}}$ staisfies the following properties
\begin{eqnarray} 
& & \alpha \wedge_{\vec{h}} \beta = 
	(-1)^{|\alpha||\beta|} \beta \wedge_{\vec{h}} \alpha,  
	\label{eqn:com} \\
& & (\alpha \wedge_{\vec{h}} \beta) \wedge_{\vec{h}} \gamma 
	= \alpha \wedge_{\vec{h}} (\beta \wedge_{\vec{h}} \gamma),
	\label{eqn:ass} 
\end{eqnarray}
for  all $\alpha, \beta, \gamma \in \Lambda_{\vec{h}}(V^*)$.
\end{theorem}

\begin{proof}
We regard $(\ref{eqn:com})$ and $(\ref{eqn:ass})$ as polynomial 
equations in $h_1, \cdots, h_m$.
For any values of  $h_1, \cdots, h_m$ in $\bk$,
set $w = h_1 w_1 + \cdots h_m w_m$. Then 
$\wedge_{\vec{h} \cdot \vec{w}} = \wedge_w$.
By Theorem \ref{thm:ass} and the remark following it,
$(\ref{eqn:com})$ and $(\ref{eqn:ass})$ hold for 
$\wedge_{\vec{h} \cdot \vec{w}}$.
Therefore, they hold as polynomial equations. 
\end{proof}

\medskip
\begin{center}
{\bf \large Part II. Geometric Applications}
\end{center}

\section{Quantum de Rham complex}

In this section, we define quantum exterior differential operator
on  Poisson manifolds. 
We follow the original route to its discovery, using
a Poisson connection, which only exists on regular Poisson manifolds. 
It has the advantage of making the verification of the desirable 
properties conceptually simple.
Then we will show that it is actually related to some well-known 
operators which can be defined on any Poisson manifolds.
Properties of the quantum exterior differential are then 
re-established using also proved properties of these operators.

Let $M$ be a smooth manifold, with a fixed bi-vector field
$w \in \Gamma(\Lambda^2(TM))$, 
then $w$ induces quantum exterior product on $\Omega(M)$
by  fiberwise quantum exterior product.
Suppose now there is a torsionless connection $\nabla$ on $TM$ 
which preserves $w$. 
Then $w$ defines a Poisson structure on $M$.
And since $w$ is parallel, it has constant rank.
Poisson manifolds whose Poisson bi-vector fields have constant ranks
are called regular Poisson manifolds. 
Conversely, any regular Poisson manifold admits a torsionless
connection which preserves the Poisson bi-vector.
Such a connection  is called a Poisson connection.
Symplectic manifolds are examples  of regular Poisson manifolds. 
See Vaisman \cite{Vai} p. 11 and p. 29 for more details.

Let $\{e_i\}$ be a local frame of $TM$ near $x \in M$,
and $\{e^i \}$ be the dual frame of $T^*M$. 
Define $d_h: \Omega(M) \rightarrow \Omega(M)[h]$ by 
$$d_h \alpha = e^i \wedge_h \nabla_{e_i} \alpha,$$
for $\alpha \in \Omega(M)$.
This definition clearly does not depend on 
the choice of the basis $\{ e_i \}$.
It follows then that we can use some particularly chosen frame 
to simplify the calculations.
Near each $x \in M$, we will use the normal coordinates 
with respect to $\nabla$.
As in Riemannian geometry, we consider the geodesics through $x$ with
respect to $\nabla$, i.e. smooth curves $c: (-1, 1) \rightarrow M$,
$c(0) = x$ and $\nabla_{c'(t)}c'(t) = 0$. 
Given any basis of $T_xM$, one can use parallel transport along 
the geodesics starting from $x$
to construct a local frame $\{ e_i \}$ near $x$. 
It then follows that $\nabla_{e_i} e_j = 0$ at $x$.
Similarly for the dual frame $\{ e^i \}$, if we still use 
$\nabla$ to denote the induced connection on $T^*M$, we have
$\nabla_{e_i} e^j = 0$ at $x$.
Since $\nabla$ is torsion free, we also have
$$[e_i, e_j] = \nabla_{e_i} e_j - \nabla_{e_j} e_i = 0, $$
at $x$. Furthermore, let $w^{jk} = w(e^j, e^k)$,
since $\nabla$ preserve $w$, we have
$$\nabla_{e_i} w^{jk} =  w(\nabla_{e_i}e^j, e^k) 
	+ w(e^j, \nabla_{e_i} e^k ),$$
so $\nabla_{e_i} w_{jk} = 0$ at $x$.
Given the torsion-free connection on $TM$, 
let $R_{X, Y} Z= \nabla_X \nabla_Y Z - \nabla_Y \nabla_X Z 
	- \nabla_{[X, Y]}$ 
be its curvature. 
It is well-known that Bianchi identity holds for torsionless
connections on the tangent bundle, i.e.,
$R_{(X, Y)}Z + R_{(Y, Z)}X + R_{(Z, X)}Y = 0$ for any vector fields
$X, Y, Z$ on $M$. 
Let $R_{(e_i, e_j)} e_k = R_{ij, k}^l e_l$,
then we have
$$R^l_{ij, k} + R^l_{jk, i} + R^l_{ki, j} = 0,$$
for any $i, j, k$  and $l$.
Denote by $\tilde{R}$ the curvature of the induced connection on $T^*M$,
let $\tilde{R}_{(e_i, e_j)}e^l = \tilde{R}_{ij, k}^l e^k$. 
Then it is routine to see that $\tilde{R}_{ij, k}^l = R_{ij, k}^l$.
Therefore, we have 
$$\tilde{R}^l_{ij, k} + \tilde{R}^l_{jk, i} + \tilde{R}^l_{ki, j} = 0,$$
for any $i, j, k$  and $l$.

\begin{theorem} \label{thm:OK}
For any $w \in \Gamma(\Lambda^2(TM))$ and a 
torsion-free connection $\nabla$ which preserves $w$, 
$d_h: \Omega(M) \rightarrow \Omega(M)[h]$ defined above 
can be extended to operators
\begin{eqnarray*}
& d_h: \Omega(M)[h] \rightarrow \Omega(M)[h], \\
& d_h: \Omega(M)[h, h^{-1}] \rightarrow \Omega(M)[h, h^{-1}]
\end{eqnarray*}
as derivations, i.e.,  
\begin{eqnarray} \label{eqn:derivation}
d_h(\alpha \wh \beta) = 
	(d_h \alpha) \wh \beta + (-1)^{|\alpha|}\alpha \wh (d_h \beta), 
	\label{derivation}
\end{eqnarray}
for $\alpha, \beta$ both in $\Omega(M)[h]$ or both in $\Omega(M)[h, h^{-1}]$.
Furthermore, $d_h^2 = 0$. 
\end{theorem}

\begin{proof} We only need to check at each $x \in M$, where we can use the
normal coordinates as above.
We write 
$\alpha = \alpha_I e^I$, $\beta = \beta_J e^J$. 
It is clear that $d_h e^I = 0$, $d_h e^J = 0$ at $x$.
This implies that at $x$, $d_h \alpha = d\alpha_I \wh e^I$,
and $d_h\beta = d\beta_J \wh e^J$.
Since $e^I \wh e^J$ is a sum of products of  $e^i$'s with $w_{ij}$'s,
$\nabla (e^I \wh e^J ) = 0$, and therefore,
$d_h(e^I \wh e^J ) = 0$ at $x$.  
By associativity of the quantum multiplication,
we have at $x$, 
\begin{eqnarray*}
&   & d_h(\alpha \wh \beta) \\
& = & d_h(\alpha_I \beta_J e^I \wh e^J) 
	= d(\alpha_I\beta_J) \wh (e^I \wh e^J) \\
& = &  d\alpha_I \wh (e^I \wh (\beta_J e^J)) + 
	\alpha_I (d\beta_J \wh e^I) \wh e^J \\
& = &  (d\alpha_I \wh e^I) \wh (\beta_J e^J) +
	(-1)^{|\alpha|} \alpha_I (e_I \wh d\beta_J) \wh e^J  \\
& = & (d_h \alpha) \wh \beta + (-1)^{|\alpha|} 
	(\alpha_I e_I) \wh (d\beta_J \wh e^J ) \\
& = & (d_h \alpha) \wh \beta + (-1)^{|\alpha|} \alpha \wh (d_h\beta).
\end{eqnarray*}
This proves $(\ref{eqn:derivation})$.
Taking $d_h$ on both sides of $(\ref{eqn:derivation})$,
one sees that 
$$d_h^2(\alpha \wh \beta) = (d_h^2\alpha) \wh \beta + 
\alpha \wh d_h^2(\beta).$$
Hence to prove  $d_h^2 = 0$,  
it suffices to verify it on $\Omega^0(M)$ and $\Omega^1(M)$.
Let $f \in \Omega^0(M)$, then $d_h f = df = (e_jf) e^j$,  and at $x$,
\begin{eqnarray*}
&   & d_h^2 f = e^i \wh \nabla_{e_i} ((e_jf) e^j)  \\
& = & e^i \wh ((e_ie_j f) e^j + (e_jf) \nabla_{e_i} e^j ) 
	=  (e_ie_jf) ( e^i \wedge e^j + h w^{ij} ) \\
& = & \sum_{i < j} [e_i, e_j]f \cdot (e^i \wedge e^j + h w^{ij}) 
	= 0.
\end{eqnarray*}
For $1$-forms, without loss of generality,
we can take $\alpha = e^l$. 
Then $d_h e^l = e^j \wh \nabla_{e_j} e^l$.
We claim that at $x$,
$$\nabla_{e_i} (e^j \wh \nabla_{e_j} e^l) 
 = e^j \wh \nabla_{e_i}\nabla_{e_j}e^l.$$
It follows from the claim that
\begin{eqnarray*}
&   & d_h^2 e^l = e^i \wh \nabla_{e_i} (e^j \wh \nabla_{e_j} e^l) \\
& = & e^i \wh (e^j \wh \nabla_{e_i}\nabla_{e_j}e^l) 
 	=  \tilde{R}^l_{ij, k} e^i \wh e^j \wh e^k \\
& = & 2 \sum_{i < j < k} (\tilde{R}^l_{ij, k} 
	+ \tilde{R}^l_{jk, i} + \tilde{R}^l_{ki, j}) e^i \wh e^j \wh e^k = 0.
\end{eqnarray*}
Now we prove the claim. 
Let $\nabla_{e_j} e^l =\Gamma_{j, k}^l e^k$, 
then at $x$, we have
\begin{eqnarray*}
& & \nabla_{e_i} (e^j \wh \nabla_{e_j} e^l) 
= \nabla_{e_i} (e^j \wedge \nabla_{e_j} e^l +
	h w^{jk} \Gamma_{j, k}^l)  \\
& = & \nabla_{e_i} e^j \wedge \nabla_{e_j} e^l
	+ e^j \wedge \nabla_{e_i} \nabla_{e_j} e^l
	+ h (\nabla_{e_i}w^{jk})\Gamma_{j, k}^l
	+ h w^{jk} (\nabla_{e_i}\Gamma_{j, k}^l) \\
& = & e^j \wedge \nabla_{e_i} \nabla_{e_j} e^l
	+ h w^{jk} (\nabla_{e_i}\Gamma_{j, k}^l) 
\end{eqnarray*}
On the other hand,  at $x$, 
$\nabla_{e_i}\nabla_{e_j} e^l 
= \nabla_{e_i}(\Gamma_{j, k}^l e^k) 
= (\nabla_{e_i}\Gamma_{j, k}^l)e^k + \Gamma_{j, k}^l \nabla_{e_i}e^k 
= (\nabla_{e_i}\Gamma_{j, k}^l)e^k$.
And so 
$$e^j \wh \nabla_{e_i} \nabla_{e_j} e^l =  
e^j \wedge \nabla_{e_i} \nabla_{e_j} e^l
	+ h w^{jk} (\nabla_{e_i}\Gamma_{j, k}^l).$$
The claim is proved. 
\end{proof}

\begin{remark}
In the above proof, we use the supercommutativity of quantum 
exterior product in an essential way.
This explains why we cannot define quantum de Rham cohomology 
using the Riemannian metric.
\end{remark}

It is instructive to compare with the classical objects in
Riemannian geometry. 
By Theorem II.5.12 and Lemma II.5.13 in 
Michelsohn-Lawson \cite{Law-Mic},
when a Riemannian metric $g$ is used, and $h = 1$, 
then $d_h$ in this context is $d - d^*$, 
where $d$ is the exterior differential,
and $d^*$ its formal adjoint, given by $d^*= -*d*$, 
where $*: \Lambda(T^*M) \rightarrow \Lambda(T^*M)$ is the Hodge
star operator defined by the Riemannian metric.
Furthermore, if $\nabla$ is the Levi-Civita connection for $g$,
\begin{eqnarray*}
d \alpha & = & e^j \wedge \nabla_{e_j} \alpha, \\
d^* \alpha & = & - \sum_j e_j \vdash \nabla_{e_j} \alpha,
\end{eqnarray*}
for local orthonormal frame $\{e_1, \cdots, e_n \}$.
In Poisson geometry, Koszul \cite{Kos} introduced a codifferential 
$$\delta: \Lambda^*(T^*M) \rightarrow \Lambda^{*-1}(T^*M)$$
for any Poisson manifold with  bi-vector field $w$,
$$\delta \alpha = w \vdash (d\alpha) - d (w \vdash \alpha).$$
He also proved  that $\delta^2 = 0$, $d \delta + \delta d = 0$. 
(Koszul used letter $\Delta$ for $\delta$.)
The complex $(\Omega^*(T^*M), \delta)$ is called the canonical 
complex of the Poisson manifold, 
its homology $PH_*(M)$ is called the Poisson homology.
When $(M^{2n}, \omega)$ is a symplectic manifold, 
Brylinski \cite{Bry} defined an operator $*_{w}:
\Omega^k(M) \rightarrow \Omega^{2n-k}(M)$,
an analogue of Hodge $*$-star operator. 
He identified $\delta$ on $\Omega^k(M)$ 
with $(-1)^{k+1}*_{w}d *_{w}$. 
A calculation similar to the Riemannian case 
(see Vaisman \cite{Vai}, Remark 1.16) yields that
\begin{eqnarray} \label{eqn:delta}
(\delta \alpha)_{i_2 \cdots  i_k} = 
	- w^{pq} \nabla_q\alpha_{pi_2\cdots i_k}.
\end{eqnarray}
Therefore, for symplectic manifolds, 
$d_h = d - h \delta$. 
Vaisman \cite{Vai} showed that 
$(\ref{eqn:delta})$ holds for regular Poisson manifolds.
In fact, let $\nabla$ be a torsionless connection, 
which is not  required to  preserve the Poisson bi-vector field. 
Then $(4.23')$ in Vaisman \cite{Vai} gives
the following tensorial expression:
$$(\delta \alpha)_{i_2\cdots i_k} 
	= w^{pq}\nabla_p\alpha_{qi_2\cdots i_k}
 	-\frac{1}{2} \sum_{s=2}^k(-1)^s 
 	\alpha_{uvi_2\cdots\hat{i}_s \cdots i_k}
 	\nabla_{i_s}w^{uv}.$$
When $\nabla$ is a torsionless connection which preserves 
the Poisson bi-vector $w$, 
 one  recovers $(\ref{eqn:delta})$. 
Therefore, we have

\begin{proposition}
On a regular Poisson manifold $(M, w)$, 
for any  Poisson
connection $\nabla$, $d_h = d - h \delta$.
Hence $d_h$ is independent of the choice of $\nabla$.
\end{proposition}

This important result suggests that we should have 
defined $d_h = d - h \delta$ for any Poisson manifold,
and proved the following stronger version of 
Theorem \ref{thm:OK}.

\begin{theorem} \label{thm:OK2}
For any Poisson manifold $(M, w)$, 
$d_h = d - h\delta$ satisfies $d_h^2 = 0$.
\begin{eqnarray} \label{derivation2}
d_h(\alpha \wh \beta) = 
	(d_h \alpha) \wh \beta + 
	(-1)^{|\alpha|}\alpha \wh (d_h \beta), 
\end{eqnarray}
for $\alpha, \beta$ both in $\Omega(M)[h]$ or 
both in $\Omega(M)[h, h^{-1}]$. 
\end{theorem}

\begin{proof}
Koszul \cite{Kos} proved that $\delta^2 = 0$ 
and $d \delta + \delta d = 0$, 
it then follows that $d_h^ 2 = 0$. 
We say that $D(a, b)$ holds if $(\ref{derivation2})$
holds for all $\alpha$, $\beta$ with $|\alpha| = a$, $|\beta| = b$.
Our strategy is first prove $D(1, b)$, then use 
induction to prove $D(a, b)$. 
Recall that 
$\delta \alpha= w \vdash d \alpha- d (w \vdash \alpha)$.

\noindent {\em Proof of $D(1, b)$}.
Let $\{ e^1, \cdots, e^n \}$ be a local frame of $TM$, 
$\beta \in \Omega^{[b]}_h(M)$,
then we have
\begin{eqnarray*}
d_h (e^i \wh \beta) 
	= d_h(e^i \wedge \beta + hw^{ij} e_j \vdash \beta) 
 =  d_h(e^i \wedge \beta) + h d_h(w^{ij} e_j \vdash \beta). 
\end{eqnarray*}
On the other hand,
\begin{eqnarray*}
&   &  d_h(e^i \wedge \beta) = d(e^i \wedge \beta) 
	- h (w \vdash d -d w \vdash) (e^i \wedge \beta) \\
& = & d e^i \wedge \beta - e^i \wedge d \beta 
	- h w \vdash(d e^i \wedge \beta - e^i \wedge d \beta) 
	+ h d [e^i \wedge (w \vdash \beta) 
	- w^{ij} e_j \vdash \beta] \\
& = & d e^i \wedge \beta - e^i \wedge d \beta  \\
&   & - h [(w \vdash d e^i) \wedge \beta
	+ de^i \wedge (w \vdash \beta) 
	+ w^{kl} (e_k \vdash de^i) \wedge (e_l \vdash \beta)]\\
&   &	+ h[e^i \wedge (w \vdash d \beta) 
	- w^{ij} e_j \vdash d\beta] \\
&   &   + h [de^i \wedge (w \vdash \beta) - 
		e^i \wedge d (w \vdash \beta)] \\
&   &	- hd [w^{ij} (e_j \vdash \beta)] \\
&  = &  [d e^i \wedge \beta 
	- h (w \vdash d e^i) \wedge \beta
	-h w^{kl} (e_k \vdash de^i) \wedge (e_l \vdash \beta)]\\
&   &	- [e^i \wedge d\beta + h w^{ij} e_j \vdash d\beta
	- h e^i \wedge (w \vdash d \beta) 
	+ h e^i \wedge d (w \vdash \beta)] \\
&   &   - h d [w^{ij} (e_j \vdash \beta)] \\
& = & ( de^i \wh \beta - h \delta e^i \wedge \beta) 
	- (e^i \wh d \beta - h  e^i \wedge \delta \beta) 
	- h d [w^{ij} (e_j \vdash \beta)] \\
& = & d_h e^i \wh \beta - e^i \wh d_h \beta 
	-h^2 [w^{ij} (e_j \vdash \delta \beta)]
	-h d[w^{ij} (e_j \vdash \beta)] \\
& = & d_h e^i \wh \beta - e^i \wh d_h \beta 
	- h d_h [w^{ij} (e_j \vdash \delta \beta)]
	- h^2 \delta [w^{ij} (e_j \vdash \beta)]
	- h^2 [w^{ij} (e_j \vdash \delta \beta)].
\end{eqnarray*}
Therefore, 
\begin{eqnarray*}
d_h (e^i \wh \beta)
= d_he^i \wh \beta - e^i \wh d_h \beta 
+ h^2 [\delta(w^{ij} e_j \vdash  \beta)
 + w^{ij} (e_j \vdash \delta \beta)].
\end{eqnarray*}
Now,
\begin{eqnarray*}
&   & \delta (w^{ij} e_j \vdash  \beta) \\
& = & w \vdash d (w^{ij} e_j \vdash \beta)
     - d [w \vdash (w^{ij} e_j \vdash \beta)] \\
& = &  w \vdash [d w^{ij} \wedge (e_j \vdash  \beta)
	+ w^{ij} d (e_j \vdash  \beta) ] 
	-d [w^{ij} (e_j \vdash  w \vdash \beta)] \\
& = & [d w^{ij} \wedge (e_j \vdash w \vdash \beta)
	- w^{kl}
	(e_k \vdash d w^{ij}) (e_l \vdash e_j \vdash  \beta)] \\
&   & + w^{ij} w \vdash d (e_j \vdash  \beta)
	- [d w^{ij} \wedge (e_j \vdash w \vdash \beta)
	+ w^{ij} d (e_j \vdash w \vdash \beta)]\\
& = & 	- w^{kl}
	(e_k \vdash d w^{ij}) (e_l \vdash e_j \vdash  \beta)] \\
&   &   + w^{ij} w \vdash d (e_j \vdash  \beta)
	- w^{ij}  d (e_j \vdash w \vdash \beta)]\\
& = & 	- w^{kl}
	(e_k \vdash d w^{ij}) (e_l \vdash e_j \vdash  \beta)] \\
&   &  + w^{ij} w \vdash (L_{e_j} \beta - e_j \vdash  d \beta)
	-w^{ij} [L_{e_j} (w \vdash \beta)
	- e_j \vdash d (w \vdash \beta)] \\
& = &  - w^{kl}
	(e_k \vdash d w^{ij}) (e_l \vdash e_j \vdash  \beta)] \\
&   & + w^{ij} w \vdash L_{e_j}\beta 
	- w^{ij} L_{e_j}(w \vdash \beta) 
	- w^{ij}e_j \vdash 
	[w \vdash d\beta - d(w \vdash \beta)] \\
& = & - w^{kl}
	(e_k \vdash d w^{ij}) (e_l \vdash e_j \vdash  \beta)] \\	
&   & + w^{ij} w \vdash L_{e_j}\beta 
	- w^{ij} L_{e_j}(w \vdash \beta) 
	- w^{ij}e_j \vdash \delta \beta.
\end{eqnarray*}
It is well-known that for $\alpha \in \Lambda^k(M)$, and smooth vector 
fields $X, Y, V_1, \cdots, V_k$, 
$$(L_X\alpha)(V_1, \cdots, V_k) 
= X\alpha(V_1, \cdots, V_k) - \sum_{j=1}^k\alpha (V_1, \cdots,
	[X, V_j], \cdots, V_k).$$
Therefore,
\begin{eqnarray*}
&   & L_X(Y \vdash \alpha)(V_1, \cdots, V_{k-1}) \\
& = & X((Y \vdash \alpha)(V_1, \cdots, V_{k-1})) 
- \sum_{j=1}^{k-1} 
(Y \vdash \alpha)(V_1, \cdots, [X, V_j], \cdots, V_{k-1}) \\
& = & X \alpha (Y, V_1, \cdots, V_{k-1})
- - \sum_{j=1}^{k-1} 
	\alpha(Y, V_1, \cdots, [X, V_j], \cdots, V_{k-1}) \\
& = & (L_X \alpha) (Y, V_1, \cdots, V_{k-1})
	+ \alpha([X, Y], V_1, \cdots, V_{k-1}) \\
& = & (Y \vdash L_X \alpha + [X, Y] \vdash \alpha)(V_1, \cdots, V_{k-1}).
\end{eqnarray*}
I.e., $L_X(Y \vdash \alpha)
= Y \vdash L_X \alpha + [X, Y] \vdash \alpha$.
Since we can assume that $[e_j, e_k] = 0$, we have
\begin{eqnarray*}
L_{e_j}(w \vdash \beta) = L_{e_j} \sum_{k < l} w^{kl} 
e_k \vdash e_l \vdash \beta)
= \sum_{k < l} e_j w^{kl} e_k \vdash e_l \vdash \beta 
+ \sum_{k < l} e_j w^{kl} e_k \vdash e_l \vdash L_{e_j} \beta.
\end{eqnarray*}
So we get
\begin{eqnarray*}
&   &\delta(w^{ij} e_j \vdash  \beta)
 	+ w^{ij} (e_j \vdash \delta \beta) \\
& = &  - w^{kl}
	(e_k \vdash d w^{ij}) (e_l \vdash e_j \vdash  \beta) 	
    + w^{ij} w \vdash L_{e_j}\beta 
	- w^{ij} L_{e_j}(w \vdash \beta) \\
& = & - \sum_{k, l} w^{kl} e_k w^{ij} 
	(e_l \vdash e_j \vdash  \beta)
	- \sum_j \sum_{k < l}w^{ij}  e_j w^{kl} 
	(e_k \vdash e_l \vdash \beta) \\
& = & - \sum_k \sum_{l < j}  w^{kl} e_k w^{ij} 
	(e_l \vdash e_j \vdash  \beta)
	- \sum_k \sum_{l > j}  w^{kl} e_k w^{ij} 
	(e_l \vdash e_j \vdash  \beta) \\
&   &	- \sum_j \sum_{k < l}w^{ij}  e_j w^{kl} 
	(e_k \vdash e_l \vdash \beta).
\end{eqnarray*}
For the first summation, change the indices by $k \mapsto j$, 
$j \mapsto l$, $l \mapsto k$;
for the  second summation, change the indices by $k \mapsto j$,
$j \mapsto k$. Then we get 
\begin{eqnarray*}
&   &\delta(w^{ij} e_j \vdash  \beta)
 	+ w^{ij} (e_j \vdash \delta \beta) \\
& = &	-\sum_j \sum_{k < l} w^{jk} e_j w^{il}
	(e_k \vdash e_l \vdash \beta)
	- \sum_j \sum_{l > k} w^{jl} e_j \vdash w^{ik} 
	(e_l \vdash e_k \vdash \beta) \\
&   &	- \sum_j \sum_{k < l}w^{ij}  e_j w^{kl} 
	(e_k \vdash e_l \vdash \beta) \\
& = & - \sum_j \sum_{k < l} (w^{kj}e_jw^{li}
	+ w^{lj}e_jw^{ik} + w^{ij}e_jw^{kl})
	(e_k \vdash e_l \vdash \beta) =0.
\end{eqnarray*}
The last equality holds because
$$w^{kj}e_jw^{li}
	+ w^{lj}e_jw^{ik} + w^{ij}e_jw^{kl} =  0,$$
which is equivalent to $w$ be  a Poisson bi-vector field 
(Vaisman \cite{Vai}, $(1.5)$).

\noindent {\em Proof of $D(a, b)$} 
This is in the same spirit of the proof of $A(a, b, c)$ in
Theorem \ref{thm:ass}.
Assume that $D(\leq a, b)$ has been proved.
Any $\alpha \in \Omega^{[a+1]}_h(M)$ can be locally 
written as  
$$\alpha =  e^i \wedge \alpha_i$$ 
for some local frame $\{e^1, \cdot, e^n \}$ and
some $\alpha_i \in \Omega^{[a]}_h(M)$.
Now for each $i$, $e^i \wedge \alpha_i = e^i \wh \alpha_i 
+ h f(e^i, \alpha_i)$, 
for some $f(e^i, \alpha_i) \in \Omega^{[a-1]}_h(M)$. 
Then we have
\begin{align*}
& d_h (\alpha \wh \beta) & \\
=& d_h [(e^i \wh \alpha_i + h f(e^i, \alpha_i)) \wh \beta] & \\
=& d_h [ e^i \wh (\alpha_i \wh \beta) + h f(e^i, \alpha_i) \wh \beta] & \\
=& d_h e^i \wh (\alpha_i \wh \beta) - e^i \wh d_h (\alpha_i \wh \beta)
	& \text{by $D(1, a+b - 1)$} \\
& + h d_h f(e^i, \alpha_i) \wh \beta + h (-1)^{|\alpha|-2}
	f(e^i, \alpha_i) \wh d_h \beta 
	& \text{by $D(a-1, b)$} \\
=& (d_h e^i \wh \alpha_i) \wh \beta - e^i \wh (d_h \alpha_i \wh \beta
 +(-1)^{|\alpha|-1} \alpha_i \wh d_h \beta)  & \\
& + h d_h f(e^i, \alpha_i) \wh \beta + h (-1)^{|\alpha|}
	f(e^i, \alpha_i) \wh d_h \beta & \\
= & (d_h e^i \wh \alpha_i - e^i \wh d_h \alpha_i + h d_h f(e^i, \alpha_i))
 \wh \beta & \\
&+ (-1)^{|\alpha|} (e^i \wh \alpha_i + hf(e^i, \alpha_i)) \wh d_h \beta & \\
=& d_h(e^i \wedge \alpha_i) \wh \beta 
 + (-1)^{|\alpha|} (e^i \wedge \alpha_i) \wh d_h \beta & \\
= & d_h \alpha \wh \beta + (-1)^{|\alpha|} \alpha \wh d_h \beta. &
\end{align*}
This completes the proof of Theorem \ref{thm:OK2}.
\end{proof}

\section{Quantum de Rham cohomology} \label{sec:dRCoh}

\begin{definition}  
For any Poisson manifold $(M, w)$, 
the {\em (polynomial) quantum de Rham cohomology} is defined
by 
$$Q_hH_{dR}^*(M) = \Ker d_h / \Img d_h,$$
for $d_h: \Omega(M)[h] \rightarrow \Omega(M)[h]$.
The {\em Laurent quantum de Rham cohomology} is
$$Q_{h, h^{-1}}H_{dR}^*(M) = \Ker d_h / \Img d_h,$$
for $d_h: \Omega(M)[h, h^{-1}] \rightarrow \Omega(M)[h, h^{-1}]$.
\end{definition}

As a consequence of Theorem \ref{thm:OK} and Theorem \ref{thm:OK2}, 
we have
\begin{theorem} The quantum de Rham cohomology 
$Q_hH_{dR}^*(M)$  of a Poison manifold  
has the following properties:
\begin{eqnarray*}
\alpha \wh \beta & = & (-1)^{|\alpha||\beta|}\beta \wh \alpha, \\
(\alpha  \wh \beta) \wh \gamma & = & \alpha \wh (\beta \wh \gamma),
\end{eqnarray*}
for $\alpha, \beta, \gamma \in Q_hH^*_{dR}(M)$.
Similar results hold for Laurent quantum de Rham cohomology.
\end{theorem}

The goal of this secton is to provide a method to compute the 
quantum de Rham cohomology, 
and to establish its relationship with the ordinary de Rham 
cohomology. 
The primary tool is the spectral sequences associated with any
double complex. 
This approach is motivated by Brylinski's results  \cite{Bry}.

The complex $(\Omega(M)[h], d_h)$ can be regarded as a double complex
$(C^{p, q}, -h\delta, d)$, where
$C^{p,q} = h^p\Omega^{q-p}(M)$, $p \geq 0$.
This is the analogue of Brylinski's double complex ${\cal C}_{..}(M)$ 
(\cite{Bry}, $\S 1.3$).
By the standard theory for double complex
(Bott-Tu \cite{Bot-Tu}, $\S 14$), 
there are two spectral sequences $E$ and $E'$ abutting to 
$H^*(\Omega[h], d_h) = Q_hH^*_{dR}(M)$,
with $E_1^{p, q} = h^pH^q(C^{p, *}, d) = h^pH^{q-p}_{dR}(M)$,
$(E_1')^{p, q} = h^pH^*(C^{*, q}, \delta) = h^pPH_{q-p}(M)$, $p \geq 0$. 

\begin{theorem}
For a Poisson manifold with odd Betti numbers all vanishing,
the spectral sequence $E$ degenerate at $E_1$, i.e. $d_r = 0$ for
all $r \geq 0$, hence
$Q_hH^*_{dR}(M)$ is a deformation quantization of $H^*_{dR}(M)$.
\end{theorem}

\begin{proof}
This is clear since $E_1^{p, q} = h^pH^{q-p}_{dR}(M)$  
is nontrivial only if $q - p $ is even. 
Now for $r \geq 1$, $d_r$ maps to $E_r^{p, q}$ t
o $E_r^{p+r, q-r+1}$,
so $d_r$ also maps block with $p + q$ 
even to a block with $p+q$ odd.
Therefore, $d_r = 0$ for any $r \geq 1$, 
since it maps all nontrivial blocks to trivial blocks. 
Therefore, 
$Q_hH^*_{dR}(M) \cong E_{\infty} = \oplus_{p, q} E_1^{p, q} 
= \oplus_{p, q} h^p H^{q - p}_{dR}(M) 
= H^*_{dR}(M) \otimes {\Bbb R}[h]$.
\end{proof}

Similarly, $(\Omega(M)[h, h^{-1}], d - h \delta)$ 
can be regarded as a double complex 
$(\widetilde{C}^{p, q}, -h\delta, d)$,
where $\widetilde{C}^{p, q} = h^p\Omega^{q-p}(M)$, 
$p, q \in {\Bbb Z}$. 
This is essentially Brylinski's double complex 
${\cal C}_{..}^{per}$, but with a different bi-grading.
We get two spectral sequences $\tilde{E}$ and  $\tilde{E}'$
abutting to $Q_{h, h^{-1}}H^*_{dR}(M)$, 
with $\tilde{E}_1^{p, q} = h^pH^{q-p}_{dR}(M)$, 
$(\tilde{E}_1')^{p, q} = h^pH^*(C^{*, q}, \delta) = h^pPH_{q-p}(M)$,
$p, q \in {\Bbb Z}$. 
The same proof yields

\begin{theorem}
For a  Poisson manifold with odd Betti numbers  all vanishing,
the spectral sequence $\tilde{E}$ degenerates at 
$\tilde{E}_1$, i.e. $d_r = 0$ for
all $r \geq 0$, hence
$Q_{h, h^{-1}}H^*_{dR}(M)$ is a Laurent 
deformation quantization of $H^*_{dR}(M)$.
\end{theorem}

Brylinski \cite{Bry} proved that on closed K\"{a}hler manifold $(M, \omega)$,
every de Rham cohomology class has a representative $\alpha$ such that
$d \alpha = 0$, $\delta \alpha = 0$.
This implies that 

\begin{Theorem} 
For a closed K\"{a}hler manifold $M$,
the spectral sequence $E$ degenerate at $E_1$, i.e. $d_r = 0$ for
all $r \geq 0$, hence
$Q_hH^*_{dR}(M)$ is a deformation quantization of $H^*_{dR}(M)$.
\end{Theorem}

We now assume  that $(M^{2n}, \omega)$ is 
a compact symplectic manifold without boundary.
Then Corollary 2.2.2 of Brylinski \cite{Bry} states that 
$PH_i(M) \cong H^{2n-i}_{dR}(M)$.
This is one of the main ingredient in Brylinski's proof  
of Theorem 2.3.1, which states that one of the spectral sequences 
for his double complex ${\cal C}_{..}^{per}$ degenerates at $E_1$.
Therefore, we have

\begin{theorem} \label{thm:degeneracy}
For any compact symplectic manifold without boundary,
the spectral sequences $\tilde{E}$ and $\tilde{E}'$ degenerate 
at $\tilde{E}_1$ and $\tilde{E}_1'$ respectively.
Hence
$Q_{h, h^{-1}}H^*_{dR}(M)$ is a Laurent 
deformation quantization of $H^*_{dR}(M)$.
\end{theorem}

\begin{proof}
The degeneracy of $\tilde{E}'$ is Brylinski's Theorem 2.3.1. 
The degeneracy of $\tilde{E}$ is by 
the following dimension counting argument.
Since $\dim \tilde{E}^{p, q}_{\infty} \geq \dim\tilde{E}^{p, q}_1$
for all $p, q \in {\Bbb Z}$,
with equalities hold for all $p, q$ 
if and only if $\tilde{E}$ degenerates at $\tilde{E}_1$,
we have
\begin{eqnarray*}
&  & \dim Q_{h, h^{-1}}H^m(M) = 
	\sum_{p+q = m} \dim \tilde{E}^{p, q}_{\infty} \\
& \geq &   \sum_{p+q = m} \dim\tilde{E}^{p, q}_1
 = \sum_{p+q = m} \dim H_{dR}^{q -p}(M),
\end{eqnarray*}
for all $m \in {\Bbb Z}$, with equalities hold for $m$ if and only if
$\tilde{E}$ degenerates at $\tilde{E}_1$.
On the other hand, since $\tilde{E}'$ degenerates at 
$\tilde{E}'_1$, we have
\begin{eqnarray*}
&  & \dim Q_{h, h^{-1}}H^m(M) = 
	\sum_{p+q = m} \dim (\tilde{E}'_{\infty})^{p, q} \\
& \geq &   \sum_{p+q = m} \dim (\tilde{E}_1')^{p, q}
 = \sum_{p+q = m} \dim PH_{q -p}(M) \\
& = & \sum_{p+q = m} \dim H_{dR}^{2n-(q -p)}(M)
= \sum_{p+q = m} \dim H_{dR}^{q -p}(M),
\end{eqnarray*}
for all $m \in {\Bbb Z}$.
\end{proof}

In fact, Brylinski's proof can be used to give 
a straightforward proof of the degeneracy of $\tilde{E}$.
It works for any closed Poisson manifold,
so we have

\begin{theorem} \label{thm:degeneracy2}
For any closed symplectic manifold,
the spectral sequences $\tilde{E}$ degenerates 
at $\tilde{E}_1$.
Hence
$Q_{h, h^{-1}}H^*_{dR}(M)$ is a Laurent 
deformation quantization of $H^*_{dR}(M)$.
\end{theorem}

\begin{remark}
Such results resemble  similar results for  the double complex 
and the associated spectral sequence appear in the Cartan model
of equivariant cohomology. 
Such spectral sequences for equivariant cohomology
appeared in Kalkman \cite{Kal}.
They were independently discovered by the  second author when he
prepared for a presentation for a course on Chern-Weil theory by 
Prof. Lawson in 1992.
For a compact symplectic manifold without boundary, Kirwan \cite{Kir}
(Proposition 5.8) 
proved that the equivariant cohomology of a Hamiltonian action by a
compact connected Lie group $G$ is a free $H^*(BG)$-module generated by
$H^*_{dR}(M)$. 
This result can be interpreted as saying the corresponding spectral
sequence for equivariant cohomology 
degenerates at $E_1$ for any compact symplectic manifold 
without boundary.
It is interesting to find a link between Kirwan's result with 
Theorem \ref{thm:degeneracy} and Theorem \ref{thm:degeneracy2}.
\end{remark}

\begin{remark}
If $M_1$, $M_2$ and $M_1 \times M_2$ all have the property that
the (Laurent) quantum de Rham cohomology is isomorphic to the
de Rham cohomology tensored with ${\Bbb R}[h]$ (${\Bbb R}[h, h^{-1}]$),
then from K\"{un}nneth formula for de Rham cohomology, 
one can deduce that $Q_hH^*_{dR}(M_1 \times M_2) 
\cong Q_hH^*_{dR}(M_1) \widehat{\otimes} Q_hH^*_{dR}(M_2)$ as
graded algebras. Similarly for the Laurent quantum de Rham cohomology.
For K\"{u}nneth formula for quantum cohomology via
pseudo-holomorphic curves, cf. Kontsevich-Manin \cite{Kon-Man2}
and Tian \cite{Tia}.
It seems plausible to develop a Leray spectral sequence for
symplectic fibration for (Laurent) quantum de Rham cohomology. 
\end{remark}

Brylinski \cite{Bry} asked the question whether every de Rham
class of a closed symplectic can be represented by a form  
$\alpha$ such that $d \alpha =0$, $\delta \alpha = 0$. 
For closed K\"{a}hler manifolds, Brylinski \cite{Bry} showed
that $\delta$ is essentially $d^*$ up to the type of the form 
it acts on.
Therefore, By Hodge theory, the answer to the above question
for closed K\"{a}hler manifolds is yes.
It has been answered negatively by Mathieu \cite{Mat} and
Yan \cite{Yan} negatively for general symplectic manifolds. 
Nevertheless,
Theorem \ref{thm:degeneracy2} implies that on a closed Poisson
manifold, any closed 
$\alpha \in \Omega^k(M)$ can be extended to a $d_h$-closed 
form $\alpha_h \in \Omega_{h, h^{-1}}^{[k]}(M)$.

\section{Quantum Hard Lefschetz Theorem}
\label{sec:Lefschetz}

The symplectic adjoint of $d_h$ is $\delta = \delta - h^{-1}d 
= h^{-1}d_h$, 
hence every quantum de Rham class is represented by a quantum
symplectic harmonic form (in the sense that $d_h \alpha =0$,
$\delta_h \alpha = 0$) for a trivial reason. 
This  aspect of quantum Hodge theory has no analogue in the
traditional approach to Hodge theory.
An important result in Hodge theory on closed K\"{a}hler manifolds
is the Hard Lefschetz theorem (Griffiths-Harris \cite{Gri-Har}, p. 122)
which states that
for a closed K\"{a}hler manifold $(M, \omega)$ of complex dimension $n$,
the map 
$$L^k: H^{n - k}(M) \rightarrow H^{n +k}(M)$$
is an isomorphism for all $k \leq n$,
where $L$ is the map given by wedge product with the K\"{a}hler form 
$\omega$.
Furthermore, if one defines the primitive cohomology 
$$P^{k-1}(M) = \Ker L^{k+1}: H^{n-K}(M) \rightarrow
	H^{n+k+2}(M),$$
then one has 
$$H^m(M) = \oplus_k L^kP^{m-2k}(M),$$
called the Lefschetz decomposition, which is compatible with the
Hodge decomposition.
This theorem is proved using results concerning finite dimensional
representations of $sl(2, {\Bbb C})$,
an idea attributed to Chern.
This no longer holds for a general symplectic manifold. 
Mathieu \cite{Mat} and Yan \cite{Yan} proved the following theorem by
different methods:

\begin{theorem}
let $(M^{2n}, \omega)$ be a symplectic manifold of dimension $2n$. 
Then the following two properties of $M$ are equivalent:
\begin{enumerate}
\item Any de Rham cohomology class of $M$ can be represented by a 
symplectic harmonic differential form.
\item For any $k \leq n$, the map $L^k: H^{n-k}(M) \rightarrow
H^{n+k}(M)$ is surjective.
\end{enumerate}
\end{theorem}

Mathieu's proof involves representation theory of 
quivers and Lie superalgebras,
Yan's proof is along the lines of the standard theory,
by considering a special class of infinite dimensional 
representations of $sl(2, {\Bbb C})$.
Motivated by all these works, 
we will study some Lie algebras of some operators
acting on the quantum exterior algebra and 
the quantum de Rham cohomology.
First notice that  $\Lambda_h^{[n-k]}(V^*)$ and 
$\Lambda_h^{[n+k]}(V^*)$ do not have the same dimension when 
$k > 0$,
so we will work with $\Lambda_{h, h^{-1}}(V^*)$.

Let $V$ be a $2n$-dimensional vector space 
over a field $\bk$ of characteristic zero.
$\omega \in \Lambda^2(V^*)$ a $\bk$-symplectic $2$-form.
Since $\omega$ is anti-symmetric, care has to be taken 
in raising or lowering the indices.
Our convention is as follows:  
$\omega$ induces an isomorphism $\sharp: V^* \rightarrow V$
by $\omega(v, \alpha^{\sharp}) = \alpha(v)$, for 
$\alpha \in V^*$, $v \in V$. 
Its inverse is denoted by $\flat: V \rightarrow V^*$.
Let $\{e_1, e_2, \cdots, e_{2n-1}, e_{2n} \}$ be a basis of $V$,
$\omega_{kl} = \omega(e_k, e_l)$. 
Let $(\omega^{kl})$ be the inverse matrix of $(\omega_{kl})$,
i.e. $\omega^{jk}w_{kl} = \delta^j_l$,
$\omega_{jk}\omega^{kl} = \delta^k_l$.
Then $(e^k)^{\sharp} = \omega^{lk}e_l$,
$e_l^{\flat} = \omega_{kl}e^k$. 
The musical isomorphism $\sharp$ induces an isomorphism
$\sharp: \Lambda^2(V^*) \rightarrow \Lambda^2(V)$ by 
$$(\phi_1 \wedge \phi_2)^{\sharp} = \phi_1^\sharp \wedge \phi_2^\sharp,$$
for $\phi_1, \phi_2 \in V^*$.
Let $w = \omega^\sharp \in \Lambda^2(V)$.
Then we have
\begin{eqnarray*}
w= \frac{1}{2} w_{kl} (e^k)^{\sharp} \wedge (e^l)^\sharp
= \frac{1}{2} w_{kl} w^{pk} w^{ql} e_p \wedge e_q 
=  \frac{1}{2} \delta^p_l w^{lq} e_p \wedge e_q 
=  \frac{1}{2} w^{pq} e_p \wedge e_q,
\end{eqnarray*}
i.e., $w^{pq} := w(e^p, e^q) = w^{pq}$.
Let $v_{w} = w^n/n!$. 
Brylinski \cite{Bry} defined the symplectic star operator
$$*: \Lambda^k(V^*) \rightarrow \Lambda^{2n-k}(V^*)$$ 
by $\beta \wedge *\alpha = \lambda^k(w)(\beta, \alpha) v_{w}$,
for all $\alpha, \beta \in \Lambda^k(V^*)$. 
He also showed that $*^2 = 0$. 
We define operators $L$,  $L^*$, $K$, $K^*$ and $A$ as follows:
\begin{eqnarray*}
L(\alpha) = w \wedge \alpha, & L^* = -*L*, \\ 
K(\alpha) = e^j \wedge (e_j \vdash \alpha), & K^* = -*K*.  
\end{eqnarray*}

\begin{lemma} \label{lemma:identities}
We have the following identities:
\begin{enumerate}
\item $L^* \alpha = w \vdash \alpha$.
\item $K(\alpha) = k \alpha$, if $\alpha \in \Lambda^k(V^*)$,
 hence $K^* = K - 2n$, $[K, K^*] = 0$.
\item $[L, K] = -2L$, $[L, K^*] = -2L$,  
$[L^*, K] = 2L^*$, $[L^*, K^*] = 2L^*$.
\item $[L, L^*] \alpha = (k - n) \alpha$, 
for $\alpha \in \Lambda^k(V^*)$. 
\end{enumerate}

\end{lemma}

\begin{proof} The first identity has been proved by Yan \cite{Yan}. 
The rest are trivial.
\end{proof}

Set $A = - \frac{1}{2}(K + K^*)$,
then we have $A(\alpha) = (n -k) \alpha$, 
for $\alpha \in \Lambda^k(V^*)$. 
Furthermore, 
\begin{eqnarray*}
[L, L^*] = A, & [L, A] =  2L, & [L^*, A] = -2L^*.
\end{eqnarray*}
This is Corollary 1.6 in Yan \cite{Yan}.
Now we define $L_h: \Lambda_{h, h^{-1}}(V) \rightarrow \Lambda_{h, h^{-1}}(V)$
by $L_h(\alpha) = w \wh \alpha$.
We extend $*$ to $\Lambda_{h, h^{-1}}(V)$ by setting $*h = h^{-1}$,
and $*h^{-1} = h$.
Define $L_h^* = -*L*$, then we have 

\begin{lemma} \label{lemma:Lh}
 We have
\begin{enumerate}
\item $L_h = L + h K + h^2 L^*$.
\item $L^*_h = L^* + h^{-1} K^* + h^{-2}L$.
\end{enumerate}
\end{lemma}

\begin{proof} Recall that $w^{ij} = w^{ij}$.
\begin{eqnarray*}  
&   & L_h(\alpha)  =  w \wh \alpha  \\
& = & w \wedge \alpha 
+ h w^{ij} (w \dashv e_i) \wedge (e_j \vdash \alpha)
+ \frac{h^2}{2!} w^{i_1j_1}w^{i_2j_2} (w \dashv e_{i_1} \dashv e_{i_2})
(e_{j_2} \vdash e_{j_1} \alpha) \\
& = & w \wedge \alpha 
	+ h w^{ij} w_{ki} e^k \wedge (e_j \vdash \alpha)
	+ \frac{h^2}{2!} w^{i_1j_1} w^{i_2j_2} 
	w_{i_2i_1}  
	(e_{j_2} \vdash e_{j_1} \vdash\alpha)\\
& = & w \wedge \alpha 
	+ h \delta^j_k e^k \wedge (e_j \vdash \alpha)
	+ \frac{h^2}{2} \delta^{j_1}_{i_2} w^{i_2j_2} 
	(e_{j_2} \vdash e_{j_1} \vdash\alpha) \\
& = & w \wedge \alpha 
	+ h e^j \wedge (e_j \vdash \alpha)
	+ \frac{h^2}{2} w^{j_1j_2} 
	(e_{j_2} \vdash e_{j_1} \vdash\alpha) \\
& = & w \wedge \alpha 
	+ h e^j \wedge (e_j \vdash \alpha)
	+ h^2 (w \vdash \alpha) \\
& = & L(\alpha) + h K(\alpha) + h^2 L^*(\alpha).
\end{eqnarray*}
The second identity follows from the first and Lemma 
\ref{lemma:identities}.
\end{proof}

We define 
$A_h: \Lambda_{h, h^{-1}}(V^*) \rightarrow 
	\Lambda_{h, h^{-1}}(V^*)$
by $A_h(\alpha) = (n - k) \alpha$, 
for $\alpha \in \Lambda^{[k]}_{h, h^{-1}}(V^*)$.

\begin{lemma} \label{lemma:Lie}
The following identities hold:
\begin{eqnarray*}
[L_h, L_h^*] = 0, & [L_h, A_h] = 2 L_h, & [L^*_h, A_h] = - 2 L^*_h.
\end{eqnarray*}
Furthermore, if we regard multiplications by $h$ and $h^{-1}$
as operators, then we have 
\begin{eqnarray*}
[h, h^{-1}] = 0, & [L_h, h^{\pm 1}] = [L_h^*, h^{\pm 1}] = 0, &
[A_h, h^{\pm}] = \pm 2 h^{\pm 1}.
\end{eqnarray*}
Therefore, if $\fg$ is the Lie algebra with three generators
$H, X, Y$, such that
\begin{eqnarray*}
[X, Y]  = 0, & [X, H] = 2X, & [Y, H] = - 2Y,
\end{eqnarray*}
then the linear map defined by 
$X \mapsto L_h$, $Y \mapsto L_h^*$, $H \mapsto A_h$
is a representation of the Lie algebra $\fg$.
Similarly, if $\fg'$ is the Lie algebra with generators
$H, X, Y, M^+, N^-$, such that
\begin{align*}
 & [X, Y]  = 0,  & &[X, H] =  2X, & &[Y, H]  = - 2Y, \\
& [M^+, M^-]  = 0,  & &[X, M^{\pm}]   = [Y, M^{\pm}]= 0, 
&&[H, M^{\pm}]   = \pm 2 M^{\pm}, 
\end{align*}
Then the linear map defined by 
$X \mapsto L_h$, $Y \mapsto L_h^*$, $H \mapsto A_h$,
$M^+ \mapsto h$, $M^- \mapsto h^{-1}$ 
is a representation of the Lie algebra $\fg'$.
\end{lemma}

\begin{proof}
>From Lemma \ref{lemma:Lh} and $K^* = K - 2n$, we see that
\begin{eqnarray} \label{eqn:Lh*}
L_h^* = h^{-2}L_h - 2nh^{-1}.
\end{eqnarray}
 Therefore $[L_h, L_h^*] = 0$.
The other  identities are trivial.
\end{proof}

It is straightforward to verify that for any constant $t$,
the $\bk$-vector space $\fg_t$ spanned by $H, X, Y$ with 
$[\cdot, \cdot]_t: \Lambda^2(\fg_t) \rightarrow \fg_t$
such that   
\begin{eqnarray*}
[X, Y]_t  = tH, & [X, H]_t = 2X, & [Y, H]_t = - 2Y,
\end{eqnarray*}
is a Lie algebra.
Over the complex field, 
it is easy to see that for any $t \neq 0$,
$(\fg_t, [\cdot, \cdot]_t)$ is isomorphic to $sl(2, {\Bbb C})$.
When $t = 0$, it gives us the Lie algebra $\fg$ in Lemma 
\ref{lemma:Lie}.
The above discussions actually suggest the following construction.
Let $\phi: sl(2, \bk) \rightarrow End(W)$ be any representation
of $sl(2, \bk)$ (over $\bk$).
As a $\bk$-vector space, $sl(2, \bk)$ is spanned by three vectors
$H, X, Y$, such that
\begin{eqnarray*}
[X, Y]  = H, & [X, H] = 2X, & [Y, H] = - 2Y.
\end{eqnarray*}
Define the following operators on $W \otimes_\bk \bk[h, h^{-1}]$:
\begin{equation} \label{operators}
\begin{split}
L_h(\pm, p) = \phi(X) \pm h \phi(H) + h^2 \phi(Y) + ph, \\
L^*_h(\pm, q) = \phi(Y) \pm h^{-1} \phi(H) + h^{-2} \phi(X) + qh,
\end{split}
\end{equation}
and $A_h(r)$ is defined to be $\phi(H) + r$ on $W$, and
$A_h(r) (h) = 2 h$, $A_h(r)(h^{-1}) = - 2h$.
Then the linear map given by $X \mapsto L_h(\pm, p)$, 
$Y \mapsto L^*_h(\pm,q)$, $H \mapsto A_h(r)$
is a representation of $\fg$.
If we also send $M^+$ to $h$, and $M^-$ to $h^{-1}$,
then we get a representation of $\fg'$. 
In particular, if $W = \Lambda(V^*)$, 
$L = \phi(X)$, $L^* = \phi(Y)$, $A = \phi(H)$, 
then 
\begin{eqnarray*}
 L_h = L_h(-, n), & L_h^* = L_h^*(-, -n), & A_h = A_h(0).
\end{eqnarray*}
To get the analogue of Hard Lefschetz Theorem,
we will  not use the representation theory for $\fg$ or
$\fg'$.  Instead, there is a simpler algebra:
let $M_n = h^{-1}L_h$, $M_n^* = hL^*_h$,
where $2n = \dim V$. 
Then $M_n^* = M_n - 2n$.
Furthermore, $M$, $M^*$ and $A_h$ commute with each other.
Since multiplications by $h$ and $h^{-1}$ are isomorphisms
which commutes with $M_n$ 
and $M_n^*$,
it suffices to examine the representation of this 
commutative Lie algebra  
on $\Omega^{[0]}_{h, h^{-1}}(V_{n}^*)$ and
$\Omega^{[1]}_{h, h^{-1}}(V_{n}^*)$.  
So we only need to find 
the eigenvalues of $M_n$ on these spaces.
We will need the following easy
lemma in linear algebra:

\begin{lemma} \label{lemma:matrices}
Let $\{ M_n \}$ be a sequence of square 
matrices with coefficient in $\bk$
obtained in the following way:
$$
M_{n+1} = \left( \begin{array}{cc} 
M_n & -I \\ I & M_n + 2I \end{array} \right)
$$
for $n \geq 1$. 
where $I$ is the identity matrix of the same size as $M_n$.

(a) For any $\lambda \in \bk$, and $n \geq 1$, we have
$$\det (M_{n+1} + \lambda I)  
=  \det [M_n + (\lambda +1)I]^2.$$
Therefore, the eigenvalues of $M_{n+1}$ can be 
obtained by adding $1$ to that of $M_n$,
with the multiplicities doubled.

(b) For any $\lambda \in \bk$, $n \geq 0$, 
$\det (M_{n+1} + \lambda I)  
= \det [M_1 + (\lambda + n) I]^{2^n}$.
Therefore, the eigenvalues of $M_{n+1}$ 
can be obtained by adding $n$ to that of $M_1$,
with $2^n$ times the multiplicities. 
If particular, if $\det (M_1 + nI) \neq 0$ for $n \geq 0$,
then $\det M_{n+1} \neq 0$.
\end{lemma}

\begin{proof}
(a) We use the standard trick of 
making one block of the matrix zero.
Notice that if $M_1$ is a $m \times m$ matrix,
then the size of $M_n$ is $m2^{n-1} \times m2^{n-1}$.
\begin{eqnarray*}
\det (M_{n+1} + \lambda I)
& = & \det \left( \begin{array}{cc} 
M_n + \lambda I & -I \\ I 
& M_n + (\lambda + 2)I \end{array} \right) \\
& = & \det \left[
\left( \begin{array}{cc} 
M_n+ \lambda I  & -I \\ 
I & M_n + (\lambda  + 2)I \end{array} \right)
\left( \begin{array}{cc} 
I & 0 \\ M_n + \lambda I & I \end{array} \right)
\right] \\
\end{eqnarray*}
\begin{eqnarray*}
& = & \det \left( \begin{array}{cc} 
0 & -I 
\\  (M_n + (\lambda + 1)I)^2 & M_n + 2 I
\end{array} \right) \\
& = & \det [M_n + (\lambda +1)I]^2.
\end{eqnarray*}
(b) An easy consequence of (a) by induction. 
\end{proof}

\begin{remark}
It is clear that similar results hold for the sequence of 
matrices defined by
$$
M_{n+1} = \left( \begin{array}{cc} 
M_n & I \\ -I & M_n - 2I \end{array} \right)
$$
for $n \geq 1$. 
For such a sequence, we have
$\det (M_{n+1} + \lambda I)  
= \det [M_1 + (\lambda - n) I]^{2^n}$,
for any $\lambda \in \bk$, $n \geq 0$.
\end{remark}

\begin{lemma} \label{lemma:eigen}
The eigenvalues of $M_1$ on 
$\Lambda^{[0]}_{h, h^{-1}}(V^*)$
are $1 \pm \frac{\sqrt{5}}{2}$,
on $\Lambda^{[1]}_{h, h^{-1}}(V^*)$, there is only one
eigenvalue $1$ with multiplicity $2$. 
For any $n > 0$, and 
any $2(n+1)$-dimensional symplectic vector space $V_{n+1}$,
the eigenvalues of the operator $M$ on both
$\Lambda^{[0]}_{h, h^{-1}}(V^*)$
and $\Lambda^{[1]}_{h, h^{-1}}(V^*)$
are $n - \frac{\sqrt{5}}{2}$, $n$ and 
$n + \frac{\sqrt{5}}{2}$.
\end{lemma}

\begin{proof} 
We will express the operator $M$ as a matrix 
in a suitable basis.
Let $\{e^1, e^2, \cdots, e^{2n+1}, e^{2n+2} \}$ 
be a symplectonormal basis of $V_{n+1}$,
let $V_n$ be the span of the first $2n$ base vectors. 
Then $\{ h^{-k} e^{i_1} \wedge \cdots \wedge e^{i_{2k}}:
k \geq 0, i_1 < \cdots < i_{2k} \}$ is a basis
of $\Lambda^{[0]}_{h, h^{-1}}(V_n^*)$,
and $\{ h^{-k} e^{i_1} \wedge \cdots \wedge e^{i_{2k+1}}:
k \geq 0, i_1 < \cdots < i_{2k+1} \}$ is a basis
of $\Lambda^{[1]}_{h, h^{-1}}(V_n^*)$.
Let $M^0_n$ and $M^1_n$ be the matrices of $M$ for $V_n$ 
in these two bases.
Now for $V_{n+1}$, 
\begin{align*}
&h^{-k} e^{i_1} \wedge \cdots \wedge e^{i_{2k}},
&& h^{-(k+1)}e^{i_1} \wedge \cdots \wedge e^{i_{2k}} \wedge
e^{2n+1} \wedge e^{2n+2}, \\
&h^{-(k+1)} e^{2n+1} \wedge 
	e^{i_1} \wedge \cdots \wedge e^{i_{2k+1}},
&& h^{-(k+1)} e^{2n+2} \wedge 
	e^{i_1} \wedge \cdots \wedge e^{i_{2k+1}}, 
\end{align*}
$k \geq 0, i_1 < \cdots < i_{2k}$, form a basis for
$\Lambda^{[0]}_{h, h^{-1}}(V_n^*)$. 
Similarly,
\begin{align*}
&h^{-k} e^{i_1} \wedge \cdots \wedge e^{i_{2k+1}},
&& h^{-(k+1)}e^{i_1} \wedge \cdots \wedge e^{i_{2k+1}} \wedge
e^{2n+1} \wedge e^{2n+2}, \\
&h^{-k} e^{2n+1} \wedge
	e^{i_1} \wedge \cdots \wedge e^{i_{2k}},
&& h^{-k} e^{2n+2} \wedge 
	e^{i_1} \wedge \cdots \wedge e^{i_{2k}}, 
\end{align*}
$k \geq 0, i_1 < \cdots < i_{2k}$, form a basis for
$\Lambda^{[1]}_{h, h^{-1}}(V_n^*)$.
Let $M^0_{n+1}$ and $M^1_{n+1}$ be 
the matrices of $M$ in these bases.
It is straightforward to verify that
\begin{eqnarray*}
M_{n+1}^0 = \left( \begin{array}{cccc}
M_n^0 & -I & 0 & 0 \\
I & M_n^0 + 2I & 0 & 0 \\
0 & 0 & M_n^1 + I & 0 \\
0 & 0 & 0 & M_n^1 + I
\end{array} \right), \\
M^1_{n+1} = \left( \begin{array}{cccc}
M^1_n & -I & 0 & 0 \\
I & M^1_n + 2I & 0 & 0 \\
0 & 0 & M^0_n + I & 0 \\
0 & 0 & 0 & M^0_n + I
\end{array} \right).
\end{eqnarray*}
In fact, for any $\alpha \in \Lambda_{h, h^{-1}}(V_n^*)$,
we have
\begin{eqnarray*}
& & M_{n+1} (\alpha) 
	= M_n (\alpha) + \alpha, \\
& & M_{n+1}(h^{-1} e^{2n+1} \wedge e^{2n+2} \wedge\alpha)
	= - \alpha 
	+  h^{-1} e^{2n+1} \wedge e^{2n+2}
		\wedge M_n (\alpha) \\
&&  \text{\hspace{2.2in}}	+ 2  h^{-1} e^{2n+1} \wedge e^{2n+2} \alpha, \\
& & M_{n+1} (e^{2n+1} \wedge \alpha) 
	= M_n (e^{2n+1} \wedge \alpha) 
	+  e^{2n+1} \wedge \alpha, \\
& & M_{n+1} (e^{2n+2}  \wedge \alpha) 
	= M_n (e^{2n+2} \wedge \alpha) 
	+ e^{2n+2} \wedge \alpha.
\end{eqnarray*}
Furthermore, we  have 
\begin{eqnarray*}
 M_1^0 =  \left( \begin{array}{cc}
 	0 & 1 \\ 1 & 2
 \end{array} \right),
&& M_1^1 =  \left( \begin{array}{cc}
	1 & 0 \\ 0 & 1
 \end{array} \right).
\end{eqnarray*}
We then inductively work out the eigenvalues of 
$M^0_{n+1}$ and $M^1_{n+1}$ with the help of 
Lemma \ref{lemma:matrices}.
The eigenvalues  of $M_1^0$ are $1 \pm \frac{\sqrt{5}}{2}$,
$M_1^1$ has eigenvalue $1$ with multiplicity $2$.
$M_2^0$ and $M_2^1$ both have eigenvalues 
$2 \pm \frac{\sqrt{5}}{2}$ and $2$.
For $n > 2$, we obtain the eigenvalues of $M_{n+1}$ by
adding $1$ to that of $M_n$.
\end{proof}

As a consequence, we have the following
algebraic version of Quantum Hard Lefschetz Theorem

\begin{theorem} \label{thm:isom}
For a symplectic vector space $V$, the operators $L_h$ and
$L_h^*$ are  isomorphisms.
Furthermore, $\Lambda_{h, h^{-1}}(V^*)$ decomposes 
into one dimensional eigen spaces of $h^{-1}L_h$ 
(or $hL_h^*$) with nonzero eigenvalues.
\end{theorem}

\begin{proof} 
Recall that $L_h = h M_n$, $L_h^* = h^{-1}M_n^*$,
and $M_n^* = M_n - 2n$. 
\end{proof}

\begin{remark}
By the same method, it is easy to find the values of
$p$ and $q$ such that
the operators $L_h(\pm, p)$ and $L^*_h(\pm, q)$ 
defined in $(\ref{operators})$ for $W = \Lambda(V^*)$
are isomorphisms.
\end{remark}

On $\Lambda_h(V^*)$, we do not have such rich structures.
It is easy to see that $\Lambda^{[n-k]}(V^*)$ and 
$\Lambda^{[n+k]}(V^*)$ do not have the same dimension when 
$k > 0$. 

Now let $(M^{2n}, \omega)$ be a 
$2n$-dimensional symplectic manifold.
Then $L_h, L_h^*, A_h$ can be defined on $\Omega_{h, h^{-1}}(M)$
by fiberwise actions.

\begin{lemma} \label{lemma:sympl}
On a symplectic manifold $(M, \omega)$, we have
\begin{eqnarray*}
[L_h, d_h] = 0, & [L_h^*, d_h] = 0, & [A_h, d_h] = -d_h.
\end{eqnarray*}
\end{lemma}

\begin{proof}
Since $\delta w 
= w \vdash d w - d(w \vdash w) = 0$,
we have $d_h w = (d - h \delta)(w) = 0$.
Therefore, for any $\alpha \in \Omega_{h, h^{-1}}(M)$,
we have
\begin{eqnarray*}
[L_h, d_h] \alpha 
= w \wh d_h \alpha - d_h (w \wh \alpha)
= - d_h w \wh \alpha = 0.
\end{eqnarray*}
The second identity follows from the first and $(\ref{eqn:Lh*})$.
The third identity is trivial.
\end{proof}

\begin{theorem}
On a symplectic manifold $(M^{2n},  \omega)$,
$Q_{h, h^{-1}}H^*_{dR}(M)$ is a representation of 
the Lie algebras $\fg$ and $\fg'$.
\end{theorem}

\begin{proof}
If $\alpha \in \Omega_{h, h^{-1}}(M)$, 
such that $d_h \alpha =0$,
then by Lemma \ref{lemma:sympl}, 
\begin{eqnarray*}
d_h (L_h \alpha) = L_h(d_h \alpha) = 0, \\
d_h (L^*_h \alpha) = L^*_h(d_h \alpha) = 0, \\
d_h (A_h \alpha) = A_h(d_h \alpha) + d_h \alpha = 0.
\end{eqnarray*}
I.e., the action of $\fg$ maps $d_h$-closed forms to 
$d_h$-closed forms.
Similarly, for any $\beta \in \Omega_{h, h^{-1}}(M)$,
\begin{eqnarray*}
 L_h(d_h \beta) = d_h (L_h \beta) = 0, \\
L^*_h(d_h \beta) = d_h (L^*_h \beta) = 0, \\
A_h(d_h \beta) = d_h (A_h \beta -  \beta) = 0.
\end{eqnarray*}
I.e., the action of $\fg$ maps $d_h$-coboundaries to 
$d_h$-coboundaries.
Therefore, the action of $\fg$ goes down to an action on 
the cohomology.
\end{proof}

As a consequence of Theorem \ref{thm:isom}, we have

\begin{theorem} (Quantum Hard Lefschetz Theorem) 
For any symplectic manifold $(M^{2n},  \omega)$,
its Laurent quantum de Rham cohomology 
$Q_{h, h^{-1}}H^*_{dR}(M)$ decomposes into 
one-dimensional eigenspaces of  the operator $h^{-1}L_h$
(or $hL_h^*$) with nonzero eigenvalues. 
In particular, $L_h$ and $L_h^*$ are isomorphisms.
\end{theorem}

\section{Quantum Dolbeault cohomology}

Let $(M, w)$ be a Poisson manifold, which admits an
almost complex structure $J$ which preserves $w$.
Assume that there is a torsionless connection $\nabla$ on $TM$,
such that $\nabla w = 0$, $\nabla J = 0$,
then $(M, w)$ is regular Poisson, and $J$ is integrable.
(As an example, consider a K\"{a}hler manifold with its Levi-Civita
connection.)
Complexify $d_h: \Omega_h(M) \rightarrow \Omega_h(M)$,
we get a decomposition 
$${\Bbb C}\Omega_h(M) = \oplus_{p, q} \Omega^{[p, q]}_h(M),$$
and correspondingly $d_h = \partial_h + \overline{\partial}_h$, 
where
\begin{eqnarray*} 
\partial_h \alpha & = & (e^i)^{1, 0} \wh 
	\nabla_{e_i^{1, 0}} \alpha, \\
\overline{\partial}_h \alpha & = & (e^i)^{0, 1} \wh 
	\nabla_{e_i^{0, 1}} \alpha,
\end{eqnarray*}
for any $\alpha \in {\Bbb C}\Omega_h(M)$.
It is clear that 
\begin{eqnarray*}
\partial_h \Omega^{[p, q]}_h(M) \subset \Omega^{[p+1, q]}_h(M), 
& \overline{\partial}_h\Omega^{[p, q]}_h(M) 
	\subset \Omega^{[p, q+1]}_h(M).
\end{eqnarray*}
Now $0 = d_h^2 = \partial_h^2 + 
	(\partial_h\overline{\partial}_h 
	+ \overline{\partial}_h\partial_h)
	+ \overline{\partial}_h^2$,
since they have bi-degrees $(2, 0)$, 
$(1, 1)$ and $(0, 2)$ respectively,
we have
\begin{eqnarray} \label{double}
\partial_h^2 = 0, & 
\partial_h\overline{\partial}_h 
	+ \overline{\partial}_h\partial_h = 0, &
\overline{\partial}_h^2 = 0.
\end{eqnarray}
Similar to $\S 2$,  the use of the connection is only an
expedient way of definition.
On a complex manifold $(M, J)$ with a Poisson structure $w$, 
not necessarily regular, such that
$J$ preserves $w$,
define $\delta^{-1, 0}: \Omega^{p, q}(M) 
	\rightarrow \Omega^{p-1, q}(M)$
and $\delta^{0, -1}(M): \Omega^{p, q}(M) 
	\rightarrow \Omega^{p, q-1}(M)$ by
\begin{eqnarray*}
\delta^{0, -1} \alpha  & = & w\vdash (\partial \alpha) 
	- \partial (w \vdash \alpha), \\
\delta^{-1, 0} \alpha  & = & 
	w\vdash (\overline{\partial} \alpha) 
	- \overline{\partial} (w \vdash \alpha),
\end{eqnarray*}
for $\alpha \in \Omega^{p, q}(M)$.
It is easy to see that for regular Poisson manifolds,
$\partial_h = \partial - h \delta^{0, -1}$,
and $\overline{\partial}_h = \overline{\partial} - h \delta^{-1, 0}$.
So we will take these as the definitions for 
$\partial_h$ and $\overline{\partial}_h$ on general Poisson manifolds.
It is clear that $(\ref{double})$ still holds.
We call 
$$Q_hH^{p, *}(M) = H(\Omega^{[p, *]}_h(M), \overline{\partial}_h)$$
the quantum Dolbeault cohomology.
We can get the quantum version of 
the usual Fr\"{o}lich spectral sequence as follows. 
>From $(\ref{double})$, we get a double complex 
$(\Omega^{[*, *]}_h(M), \partial_h, \overline{\partial}_h)$, 
whose associated complex is $({\Bbb C}\Omega^*_h(M), d_h)$, 
one of the standard spectral sequences  ${\Bbb C}Q_hH(M)$ has 
$E_1^{p, q} = Q_hH^{p, q}(M)$.
Now $0 = \overline{\partial}_h^2 = \overline{\partial}^2 
- h(\overline{\partial} \delta^{-1, 0} 
	+ \delta^{-1, 0} \overline{\partial}) 
+ h^2 (\delta^{-1, 0})^2$, so
\begin{eqnarray*}
 \overline{\partial}^2 = 0, &
 \overline{\partial} \delta^{-1, 0} + \delta^{-1, 0} \overline{\partial} = 0,
 & (\delta^{-1, 0})^2 = 0.
\end{eqnarray*}
So we get a double complex $(C^{p, q} = h^m\Omega^{p-m, n}(M), 
\overline{\partial},
-h\delta^{-1, 0})$.
It has two associated spectral sequences abutting to $Q_hH^{p, *}(M)$.
Taking cohomology in $\overline{\partial}$ first,
we get a spectral sequence with $E_1^{mn} = h^mH^{p-m, n}(M)$. 
Similar to Theorem 3.2 and Theorem 3.3, we get

\begin{theorem}
When $M$ is a compact K\"{a}hler manifold with  vanishing odd Betti
numbers, the spectral sequence of 
$(C^{m, n} = h^m\Omega^{p-m, n}(M), 
\overline{\partial},
-h\delta^{-1, 0})$ with  $E_1^{mn} = h^mH^{p-m, n}(M)$ degenerate at $E_1$. 
Hence $Q_hH^{p, q}(M) = \oplus_{k \geq 0} h^kH^{p-k, q-k}(M)$.
\end{theorem}

One can also define Laurent quantum Dolbeault cohomology 
$Q_{h, h^{-1}}H^{p, q}(M)$, and consider the corresponding spectral 
sequences.

\begin{theorem} 
For a closed K\"{a}hler manifold $M$, we have
\begin{eqnarray*}
& {\Bbb C}Q_{h, h^{-1}}H^n_{dR}(M) 
	= \oplus_{p + q = n} Q_{h, h^{-1}}H^{p, q}(M), \\
& {\Bbb C}Q_{h, h^{-1}}H^{p, q}(M) = \oplus_{k \in {\Bbb Z}} 
	h^k Q_{h, h^{-1}}H^{p-k, q-k}(M).
\end{eqnarray*}
\end{theorem}

\begin{proof}
By Theorem \ref{thm:degeneracy}, 
$$Q_{h, h^{-1}}H^n_{dR}(M) 
= \oplus_{k \in {\Bbb Z}}h^k{\Bbb C}H^{n-2k}_{dR}(M).$$
By Hodge theorem, 
${\Bbb C}H^{n-2k}_{dR}(M) \cong \oplus_{p+q = n} H^{p-k, q-k}(M)$.
So we have 
\begin{eqnarray} \label{dimI}
\dim Q_{h, h^{-1}}H^n_{dR}(M) = \sum_{k \in {\Bbb Z}} 
	\sum_{p+q = n} \dim  H^{p-k, q-k}(M),
\end{eqnarray}
where all the dimensions are dimensions as complex vector spaces,
It is a sum of finitely many finite numbers.
Now there is a spectral sequence abutting to $Q_{h, h^{-1}}H^n_{dR}(M)$
with $E_1^{p, q} = Q_{h, h^{-1}}H^{p, q}(M)$.
So we have
\begin{eqnarray} \label{dimII}
\dim Q_{h, h^{-1}}H^n_{dR}(M) \leq \sum_{p+q = n} 
	\dim Q_{h, h^{-1}}H^{p, q}(M),
\end{eqnarray}
equality holds iff the spectral sequence degenerates at $E_1$.
Similarly, 
there is a spectral sequence abutting to $Q_{h, h^{-1}}H^{p, *}(M)$
with $\tilde{E}_1^{k, l} =  h^kH^{p-k, l}(M)$. 
Therefore, 
\begin{eqnarray} \label{dimIII}
\dim Q_{h, h^{-1}}H^{p, q}(M) \leq \sum_{k+l = q} \dim H^{p-k, l}(M) 
= \sum_{k \in {\Bbb Z}} \dim H^{p-k, q-k}(M).
\end{eqnarray}
Equality holds iff the spectral sequence degenerates at $\tilde{E}_1$.
Combining $(\ref{dimII})$ with $(\ref{dimIII})$, one gets
$$\dim Q_{h, h^{-1}}H^n_{dR}(M) \leq \sum_{k \in {\Bbb Z}} 
	\sum_{p+q = n} \dim  H^{p-k, q-k}(M),$$
with equality iff both $E$ and $\tilde{E}$ degenerate at $E_1$. 
Comparing with $(\ref{dimI})$, one sees that all the relevant 
spectral sequences degenerate at $E_1$.
This completes the proof.
\end{proof}

It is easy to see that the analogue of quantum Hard Lefschetz 
Theorem holds for quantum Dolbeault cohomology..

\section{Quantum integral and quantum Stokes Theorem}

Let $(M, \omega)$ be a closed $2n$-dimensional symplectic manifold.
Define an integral $\int_h: \Omega_h(M) \rightarrow {\Bbb R}[h]$
as follows.
For any $\alpha \in \Omega^{j}(M)$,
if $j$ is odd, set $\int_h \alpha = 0$;
if $j = 2n - 2k $ for some integer $k$,
set 
$$\int_h \alpha = \int_M \alpha \wedge \frac{\omega^k}{k!}.$$
Extend $\int_h$ to $\Omega_h(M)$ as a ${\Bbb R}[h]$-module map.
We call $\int_h$ the {\em quantum integral}.
Straightforward calculations yield the following

\begin{lemma}
For $\alpha, \beta \in \Omega(M)$, we have 
$$w \vdash (\alpha \wedge \beta) = (w \vdash \alpha) \wedge \beta
+ 2 w^{ij} (e_i \vdash \alpha) \wedge ( e_j \vdash \beta) 
+ \alpha \wedge (w \vdash \beta).
$$
\end{lemma}

\begin{lemma} \label{lem:aux}
(i) We have
$$ w \vdash \frac{\omega^{k+1}}{(k+1)!}
= (n + k) \frac{\omega^k}{k!}.$$
(ii) For $\beta \in \Omega^p(M)$, we have
$$w^{ij} (e_i \vdash \beta) \wedge 
	(e_j \vdash \frac{\omega^{k+1}}{(k+1)!})
= (-1)^{p-1} p \beta \wedge \frac{\omega^k}{k!}.$$
\end{lemma}

\begin{theorem} (Quantum Staokes Theorem) 
 We have $\int_h d \alpha = 0$, 
	$\int_h h\delta \alpha = 0$, and therefore
$$\int_h d_h \alpha = 0.$$
\end{theorem}

\begin{proof}
We can assume that $\alpha$ has odd degree,
write
$\alpha = \sum_{k=0}^{n} \frac{h^k}{k!} \alpha_{2n-1- 2k}$,
where $\deg (\alpha_{2n-1-2k}) = 2n - 2k -1$. 
Then we have
\begin{eqnarray*}
\int_h d \alpha & = & \int_M \sum_{k=0}^{n} 
	d\alpha_{2n-1- 2k} \wedge \frac{\omega^k}{k!}\\
& = & \int_M \sum_{k=0}^{n} d 
	(\alpha_{2n-1- 2k}\wedge \frac{\omega^k}{k!}) = 0.
\end{eqnarray*}
Recall that $\delta \alpha  = w \vdash d \alpha - d (w \vdash \alpha)$,
therefore,
\begin{eqnarray*}
&   & \int_h h\delta \alpha = \int_h w \vdash d \alpha 
  =  \int_M \sum_{k=0}^{n}  (w \vdash d \alpha_{2n-1- 2k})
	\wedge \frac{\omega^{k+1}}{(k+1)!} \\
& = & \sum_{k=0}^{n}  \int_M 
	 w \vdash (d \alpha_{2n-1- 2k} \wedge \frac{\omega^{k+1}}{(k+1)!})  \\
& - &  \sum_{k=0}^{n} \int_M 
	d \alpha_{2n-1- 2k} \wedge 
	(w \vdash\frac{\omega^{k+1}}{(k+1)!}) \\
& - &	2 \sum_{k=0}^{n} \int_M
	w^{ij} (e_i \vdash d\alpha_{2n-1- 2k}) \wedge
	(e_j \vdash \frac{\omega^{k+1}}{(k+1)!}) .
\end{eqnarray*}
The first term vanishes since 
$d \alpha_{2n-1- 2k} \wedge \omega^{k+1}$
has degree $2n + 2 > \dim (M)$. 
By Lemma \ref{lem:aux}, 
\begin{eqnarray*}
&  & \int_h h\delta \alpha \\ & = &  
	- \sum_{k=0}^{n} (n+k) \int_M 
	d \alpha_{2n-1- 2k} \wedge \frac{\omega^k}{k!}
+ 2\sum_{k=0}^{n} (2n - 2k) \int_M 
	d \alpha_{2n-1- 2k} \wedge \frac{\omega^k}{k!} = 0.
\end{eqnarray*}

\end{proof}

\section{Quantum Chern-Weil theory}

The classical constructions in Chern-Weil theory of representing
characteristic classes of a vector bundle over a smooth manifold
by  curvature expressions can be generalized in the context of
quantum de Rham cohomology. As usual, the case of a complex line
bundle is very simple. 
We will go over it first to illustrated the idea. 
Let  $L$ be  a complex line bundle on a Poisson manifold $M$.
Given a fine open  covering $\{ U_{\alpha} \}$ of $M$, i.e. each 
$U_{\alpha}$ and each $U_{\alpha} \cap U_{\beta}$ are contractible. 
Then $L|_{U_{\alpha}}$ can be trivialized by a nonvanishing
section $s_{\alpha}$, and on $U_{\alpha} \cap U_{\beta}$,
there is a smooth complex valued function $f_{\alpha\beta}$,
such that $s_{\alpha} = \exp f_{\alpha\beta} s_{\beta}$. 
Suppose that $L$ has a connection $\nabla_L$.
Then there is a complex valued 
$1$-form $\theta_{\alpha}$ on each $U_{\alpha}$, 
such that $\nabla_Ls_{\alpha} = \theta^{\alpha} \otimes s_{\alpha}$.
Then we have $\theta^{\alpha} = \theta^{\beta} +  d f_{\alpha}$. 
We defined the quantum curvature $\Omega_h$ of $L$ by: 
on each $U_{\alpha}$, $\Omega_h = d_h \theta^{\alpha}$.
Since on $U_{\alpha} \cap U_{\beta}$,
$$d_h \theta^{\alpha} = d_h \theta^{\beta} +  d_hd f_{\alpha} 
= d_h \theta^{\beta}.$$
we have $\Omega_h \in \Omega(M)[h]$.
Similarly, one can show that  
$\Omega_h$ does not depend on the choice
of the local trivializations. 
Now clearly $d_h \Omega_h = 0$, 
we call 
$$c_1(L)_h = \frac{\sqrt{-1}}{2\pi} [\Omega_h] \in Q_hH_{dR}(M)$$ 
the quantum first Chern class of $M$. 
Use the quantum multiplications, one can also define the quantum
Chern character 
$$ch(L)_h = \exp_h \left(\frac{\sqrt{-1}}{2\pi} [\Omega_h]\right)
= \sum_{n = 0}^{\infty} \frac{(\sqrt{-1})^n}{n!(2\pi)^n} 
	[ (\Omega_h)^n_h] , $$ 
where by definition, $(\alpha)^n_h = \alpha \wh \cdots \wh \alpha$
($n$-times),
$\exp_h (\alpha) = \sum_{n =)}^{\infty}
\frac{1}{n!} (\alpha)_h^n$, for $\alpha \in \Omega_h(M)$.

Now let $E \rightarrow M$ be a vector bundle over a 
Poisson manifold $M$.  
A connection on $E$ is a linear operator
$\nabla_E: \Omega^0(E) \rightarrow \Omega^1(E)$, such that
$$\nabla (\sigma \cdot f) 
= (\nabla \sigma) \cdot f + \sigma \otimes df,$$
where $\sigma$ is any section of $E$, 
and $f$ is any smooth function on $M$.
Let $\Omega^*(E)$ be the space of exterior forms with values in $E$.
We give $\Omega_h^*(E): =\Omega^*(E)[h]$ 
a structure of right $\Omega_h^*(M)$-module.
Given a connection on $E$, define the quantum covariant derivative 
$d_h^{\nabla_E}: \Omega_h^*(E) \rightarrow \Omega_h^*(E)$ as follows.
Let $\bs$ be a local frame of $E$, let $\theta$ be the 
connection $1$-form in this frame: 
$\nabla \bs = \bs \otimes \theta$, i.e.,
$$\nabla \bs_j = \sum_{k=1}^{n} \bs_k \otimes \theta_j^k.$$
If $\alpha =\bs \otimes \phi$, for some vector valued form $\phi$,
define
$$d_h^{\nabla_E} \alpha = \bs \otimes (\theta \wh \phi + d_h \phi) 
	= \sum  \bs_k \otimes ( \theta_j^k \wh \phi^j + d_h \phi^k).$$

\begin{lemma} \label{lm:indep}
The definition of $d_h$ is independent of 
the choice of the local frames.
\end{lemma}

\begin{proof}
If $\bs'$ is another local frame such that  $\bs' = \bs \cdot G$,
$\alpha = \bs' \otimes \phi'$, 
and $\nabla \bs' = \bs' \otimes \theta'$.  
Then $\phi' = G^{-1}\phi$, $\theta' = G^{-1}\theta G + G^{-1}dG$.
Hence,
\begin{eqnarray*}
&   & \bs' \otimes (\theta' \wh \phi' + d_h \phi') \\
& = & \bs G \otimes (G^{-1}\theta G \wh G^{-1}\phi 
	+ G^{-1}dG G^{-1} \wh \phi
	+ d_h (G^{-1}\phi)) \\
& = & \bs \otimes [(\theta + dG G^{-1})\wh \phi + G^{-1}dG \wh \phi 
	+ d_h\phi] \\
& = & \bs \otimes (\theta \wh \phi + d_h \phi).
\end{eqnarray*}
\end{proof}

Alternatively, let $\{e_i\}$ be a local frame of $TM$ near $x \in M$,
and $\{e^i \}$ be the dual frame of $T^*M$. Then
$$d_h^{\nabla_E} \alpha = e^i \wedge_h \nabla'_{e_i} \alpha,$$
where $\nabla'$ is the connection on $\Lambda(T^*M) \otimes E$ induced
by the admissible connection $\nabla$ on $TM$ and $\nabla_E$ on $E$.
>From the definition and Theorem \ref{thm:OK}, it is routine to 
verify the following

\begin{lemma}
The quantum covariant derivative is a $\Omega_h^*(M)$-module derivation
of degree $1$, i.e.,
$$d_h^{\nabla_E} (\Phi \wh \alpha) = (d_h^{\nabla_E} \Phi) \wh \alpha
+ (-1)^{\deg \Phi} \Phi \wh (d_h \alpha),$$
where $\Phi \in \Omega^k_h(E)$, $\alpha \in \Omega^*_h(M)$.
\end{lemma}

Notice that $\Omega_h^*(E)$ is also a right 
$\Omega^*_h(\End(E))$-module.

\begin{theorem/definition}
There is an element $R^E_h \in
\Omega^2_h(\End(E))$, such that
for each $k \geq 0$, $(d_h^{\nabla_E})^2$ on $\Omega^k_h(M)$ is
given by $(d_h^{\nabla_E})^2 \Phi = \Phi \wh R_h^E$, for
any $\Phi \in \Omega^*_h(E)$.
$R^E_h$ is called the quantum curvature of $\nabla^E$.
\end{theorem/definition}

\begin{proof}
We use the local frame $\bs$ and local connection $1$-form 
$\theta$ as above.
Then by Theorem \ref{thm:OK},
\begin{eqnarray*}
&   & (d_h^{\nabla_E})^2 \Phi = d_h^{\nabla_E} 
( \bs \otimes (\theta \wh \phi + d_h \phi)) \\
& = & \bs \otimes \{\theta \wh (\theta \wh \phi + d_h \phi)
	+ d_h (\theta \wh \phi + d_h \phi)\} \\
& = & \bs \otimes (\theta \wh \theta_h \wh \phi + \theta \wh d\phi 
	+ d_h\theta \wh \phi - \theta \wh d_h\phi )\\
& = & \bs \otimes \{ (d_h \theta + \theta \wh \theta) \wh \phi \}.
\end{eqnarray*}
For a different local frame $\bs' = \bs G$ with 
$\nabla_E \bs' = \bs' \otimes \theta'$.
A calculation as in the ordinary case shows that
$$d_h \theta' + \theta' \wh \theta' = 
G^{-1} (d_h \theta + \theta \wh \theta)G.$$
This shows that $\Theta^{\bs}_h := d_h \theta + \theta \wh \theta$ 
in different frames patches up to give us an element $R^E$ 
in $\Omega_h^2(\End(E))$.
\end{proof}

For $(n \times n)$-matrix valued differential forms 
$\alpha = (\alpha_{ij})$
and $\beta = (\beta_{ij})$, define
$$[\alpha \wh \beta]_{ij} = \sum_k (\alpha_{ik} \wh \beta_{kj} -
	\beta_{ik} \wh \alpha_{kj}).$$ 
In a local frame $\bs$, we have 
\begin{eqnarray*}
& & d^{\nabla_E}_h (d^{\nabla_E}_h)^2 \bs 
	= d^{\nabla_E}_h (\bs \otimes \Theta_h^{\bs}) 
  =  \bs \otimes (d_h\Theta^{\bs}_h + \theta \wh \Theta_h), \\
& &  (d^{\nabla_E}_h)^2 d^{\nabla_E}_h \bs  
	= (d^{\nabla_E}_h)^2 (\bs \otimes \theta) 
  =  \bs (\Theta_h \wh \theta).
\end{eqnarray*}
Since $d^{\nabla_E}_h (d^{\nabla_E}_h)^2 = 
(d^{\nabla_E}_h)^2 d^{\nabla_E}_h$,
we get
\begin{eqnarray} \label{eqn:Bianchi}
d_h \Theta_h = [\Theta_h \wedge_h \theta].
\end{eqnarray}
If $p$ is a polynomial on the space of $n\times n$-matrices, 
such that $p(G^{-1}AG) = p(A)$, 
for any invertible $n \times n$-matrix $G$,
then $p(\Theta^{\bs})$ for different frames patch up to 
a well-defined element $p(R^E) \in \Omega^*(M)[h]$.
Similar to the ordinary Chern-Weil theory, 
it is easy to see that $d_h p(R^E) = 0$.
So it defines a class in $Q_hH^*_{dR}(M)$.
The usual  construction of transgression operator  
carries over to show that
this class is independent of the choice of the connection $\nabla^E$.
In this way,  one can define quantum Chern classes, 
quantum Euler class, etc.
We will call them quantum characteristic classes.
It is clear that we can repeat the same story in Laurent case. 
Notice that in Atiyah-Singer index theorem,
the index of an elliptic operator on a closed manifold
is expressed as the integral of a power series of the 
curvature.
If we use quantum curvature and quantum exterior product in
the power series,
we then get a power series in $h$, whose $0$-th order term 
yields the ordinary index.

\subsection{Quantum equivariant de Rham cohomology}

Let $(M, w)$ be a Poisson manifold, which 
admits an action by a compact connected Lie group $G$,
such that the $G$-action preserves the Poisson bi-vector field $w$.
Let $\fg$ be the Lie algebra of $G$, $\{ \xi_a \}$ a basis of $\fg$, 
denote by $\iota_a$ the contraction by the vector field 
generated by the one parameter group corresponding to $\xi_a$,
and $L_a$ the Lie derivative by the same vector field.
Imitating the Cartan model for equivariant cohomology, 
we consider the operator $D_{hG} = d_h + \Theta^a \iota_a 
= d - h \delta + \Theta^a \iota_a$
acting on $(S(\fg^*) \otimes \Omega(M))^G[h]$.  
It is well-known that $d + \Theta^a \iota_a$ maps 
$(S(\fg^*) \otimes \Omega(M))^G$ to itself.
Since the $G$-action preserves $w$, it is easy to check that
$\delta$ also preserves $(S(\fg^*) \otimes \Omega(M))^G$.
Therefore, $D_{hG}$ is an operator from 
$(S(\fg^*) \otimes \Omega(M))^G[h]$
to itself. 
Now on $(S(\fg^*) \otimes \Omega(M))^G[h]$, we have 
\begin{eqnarray*}
D_{hG}^2 
& = & d_h^2 + (\Theta^a\iota_a)^2 + \Theta^a (d \iota_a + \iota_a d)
	- h \Theta^a (\delta \iota_a + \iota_a \delta) \\
& = &	- h \Theta^a (\delta \iota_a + \iota_a \delta).
\end{eqnarray*}
Since $\delta = \iota_{w} d - d \iota_{w}$, we have
\begin{eqnarray*}
\delta \iota_a + \iota_a \delta 
& = & \iota_{w}d \iota_a - d \iota_{w} \iota_a 
	+ \iota_a \iota_{w} d - \iota_a d \iota_{w} \\
& = & \iota_{w}d \iota_a - d \iota_a \iota_{w}
	+ \iota_{w} \iota_a d - \iota_a d \iota_{w} \\
& = & \iota_{w}L_a - L_a \iota_{w} = -\iota_{L_aw} = 0.
\end{eqnarray*}
Hence, $D_{hG}^2 =0$.
We call the cohomology 
$$Q_hH^*_{G}(M) := H^*((S(\fg^*) \otimes \Omega(M))^G[h], D_{hG})$$
the quantum equivariant de Rham cohomology. 
Similar definitions can be made using Laurent deformation. 
We will study quantum equivariant de Rham cohomology in a forthcoming
paper.

\section{Computations for some examples}

The quantum Chern-Weil theory and Theorem \ref{thm:degeneracy}
provide us with tools to compute 
the quantum de Rham cohomology rings of some important examples
of symplectic manifolds such as projective spaces, complex
Grassmannians and flag manifolds.

\begin{example} (Complex projective space) 
For any symplectic form on ${\Bbb C}{\Bbb P}_n$,
$H^*_{dr}({\Bbb C}{\Bbb P}_n)$ is the ring ${\Bbb R}[\omega]/
(\omega^{n+1} = 0)$.
By Theorem \ref{thm:degeneracy}, 
$Q_hH_{dr}^*({\Bbb C}{\Bbb P}_n) = H^*_{dr}({\Bbb C}{\Bbb P}_n)
\otimes {\Bbb R}[h]$. 
So we need to compute $\omega^k \wh \omega^l$. 
It is clear from the definition that it is 
 a linear combination of 
 $\omega^{k+l}, \omega, \cdots, \omega^{|k-l|}$, 
with coefficients polynomials of $h$.
This can be done inductively as follows:
we first compute $\omega \wh \omega^k$,
then by induction compute $(\omega)^k_h =
\omega \wh \cdots \wh \omega$ ($k$ times),
then for $k \geq 2$, reduce the computation of 
$\omega^k \wh \omega^l$  to the first compuations. 
In fact, by Darboux Theorem,  locally we write 
$\omega = e^1 \wedge e^2 + \cdots + e^{2n-1} \wedge e^{2n}$.
By results in $\S \ref{sec:Lefschetz}$, we have
\begin{eqnarray} \label{eqn:1xk}
\omega \wh \omega^k= \omega^{k+1} + 2k h \omega^k - kn h^2
\omega^{k-1},
\end{eqnarray}
for $k \geq 1$.
>From this, we inductively compute 
$(\omega)^k_h$:
if $(\omega)^k_h = a^{(k)}_k(h) \omega^k + \cdots + a^{(k)}_0(h)$ 
for some polynomials 
$a^{(k)}_k(h), \cdots, a^{(k)}_0(h) \in {\Bbb R}[h]$,
then 
$$(\omega)^{k+1}_h = \omega \wh (\omega)^k_h 
= a^{(k)}_k(h) \omega \wh \omega^k + \cdots + a^{(k)}_0(h) \omega.$$
Use $(\ref{eqn:1xk})$, $(\omega)^{k+1}_h$ can be written as a linear
combination of $\omega^l$'s. 
Now we show how to inductively reduce the calculation of 
$\omega^{k+1} \wh \omega^l$ as follows.
By $(\ref{eqn:1xk})$,
$\omega^{k+1} = \omega \wh \omega^k 
- 2kh \omega^k + knh^2\omega^{k-1}$,
then
\begin{eqnarray*}
\omega^{k+1} \wh \omega^l 
= \omega \wh (\omega^k \wh \omega^l) - 2kh \omega^k \wh \omega^l
+ knh^2\omega^{k-1} \wh \omega^l.
\end{eqnarray*}
This can be written as a linear combination of 
$\omega^{k+1+l}, \cdots, \omega^{|k+1-l|}$ by induction hypothesis
and $(\ref{eqn:1xk})$.
The result of such recursive procedures is very complicated.
But it gives us the full information about how the ring structure
is deformed.
On the other hand, there is  simple ways to give  presentations
of the deformation quantization.
Let $\omega_h = e^1 \wh e^2 + \cdots e^{2n-1} \wh e^{2n} 
= \omega -nh$.
Since $\wh$ is super-commutative,  if we 
expand $(\omega_h)^{n+1}_h$ as a sum of terms of the form
$(e^{2i_1-} \wh e^{2i_1}) 
\wh \cdots \wh (e^{2i_{n+1}-1} \wh e^{2i_{n+1}})$,
it is clear that 
\begin{eqnarray} \label{eqn:power}
(\omega_h)^{n+1}_h = 0.
\end{eqnarray} 
This does not imply that we get a trivial deformation,
since $(\cdot)^{n+1}_h$ is given by the deformed multiplication.
Equivalently,
\begin{eqnarray*}
(\omega)^{n+1}_h = \sum_{k=0}^{n}
(-1)^{n-k} h^{n+1-k} \left( 
\begin{array}{cc} n+1 \\ k \end{array} \right)
 (\omega)^{k}_h.
\end{eqnarray*}

Let $\nu$ be the tautological line bundle over ${\Bbb C}{\Bbb P}_n$,
and $Q = \underline{\Bbb C}^{n+1}/\nu$, 
where $\underline{\Bbb C}^{n+1}$ is the trivial rank $n+1$ bundle.
Then from the exact sequence  
$0 \rightarrow \nu \rightarrow  \underline{\Bbb C}^{n+1}
\rightarrow Q \rightarrow 0$,
we get $c(\nu)_h \wh c(Q)_h = c(\underline{\Bbb C}^{n+1})_h = 1$.
Therefore,
$$c(Q)_h = 1/c(\nu)_h = \sum_{j \geq 0} (-c_1(\nu)_h)_h^j.$$
Without loss of generality, we can assume that 
$[\omega]$ is $c_1(\nu)$.
It is easy to check that $c_1(\nu)_h  = -(\omega - \lambda h)$
for some constant $\lambda$.
Since $Q$ has rank $n$, 
$(\omega + \lambda h)^{n+1}_h = c_{n+1}(Q)_h = 0$.
It follows from $(\ref{eqn:power})$ that $\lambda = n$.
On the ohter hand, if $\omega$ is the K\"{a}hler form for 
Fubuni-Study metric, 
then it is possible to check that $c_1(\nu)_h =- (\omega - n h)$.
This then yields $(\ref{eqn:power})$.
 
 This example illustrates the complexity in the calculation  
of quantum multiplications in quantum de Rham cohomology.
\end{example}

\begin{example} (Complex Grassmannian) 
The same method can be used for complex Grassmannian 
$G_{r, n}({\Bbb C})$.
Let $\nu$ be the tautological vector bundle,
and $Q = \underline{\Bbb C}^n/\nu$.
Let $c_j = c_j(\nu)$, and $s_j = c_j(Q)$.
Then from the exact sequence 
$0 \rightarrow \nu \rightarrow  \underline{\Bbb C}^{n}
\rightarrow Q \rightarrow 0$, we get
$c(\nu) \wedge c(Q) = c(\underline{\Bbb C}) = 1$, i.e
\begin{eqnarray} \label{eqn:Segre}
s_j = - s_{j-1}c_1 - \cdots - s_1c_{j-1} - c_j,
\end{eqnarray}
for $j \geq 1$.
Since $Q$ has rank $n-r$, we must have
$s_j = 0$ for $j > n -r$.
In fact, 
the de Rham cohomology ring of $G_{r, n}({\Bbb C})$ is given by
(e.g. Fulton \cite{Ful}, Ex. 14.6.6)
$${\Bbb R}[c_1, \cdots, c_r]/(s_{n-r+1}, \cdots, s_n),$$
where $s_j$'s are given by $(\ref{eqn:Segre})$.
Since $H^2_{dR}(G_{k, n}({\Bbb C})$ is one-dimensional, 
given any symplectic structure $\omega$ on $G_{k, n}({\Bbb C}))$,
we may assume without loss of generality that $[\omega] = - c_1$.
Let  $c_{j, h} = c_j(\nu)_h$, $s_{j, h} = c_j(Q)_h$.
Then using quantum Chern classes, we get $s_{j, h} = 0$,
for $j = n - r + 1, \cdots, n$,
where $s_j$'s are given by 
\begin{equation}
s_{j, h} = - s_{j-1, h} \wh c_{1, h} - \cdots 
	- s_{1, h} \wh c_{j-1, h} - c_{j, h}. \tag{\ref{eqn:Segre}'}
\end{equation}
Therefore, we need to compute $c_{j, h}$ and their multiplications.
\end{example}

\begin{example} (Complex flag manifold)
Let $F_{n+1}$denote the manifold of complete flags 
$0 = V_0 \subset V_1 \subset \cdots \subset V_{n-1} \subset V_n = {\Bbb C}^n$,
where each $V_j$ is a subspace of dimension 
$j$, for $j = 1, \cdots, n$.
There are tautological line bundles $L_j$ on $F_{n+1}$, whose 
fiber at a flag $V_1 \subset \cdots \subset V_{n-1} \subset V_n$
is $V_j/V_{j-1}$.
Then it is clear that
$$L_1 \oplus \cdots \oplus L_n = \underline{\Bbb C}^n.$$
This is a special case of splitting principle 
(see e.g. Bott-Tu \cite{Bot-Tu}, $\S 21$, especially
the proof of Proposition 21.15), which states
that for any complex vector bundle $E$ of rank $n$, 
if $\pi: F(E) \rightarrow M$ is 
the flag bundle associated with $E$,
then $\pi^*E = L_1 \oplus \cdots \oplus L_n$,
where $L_j = V_j(E)/V_{j-1}(E)$.
Therefore, we have
\begin{eqnarray} \label{eqn:flag}
c(L_1) \cdots c(L_n) = 1.
\end{eqnarray}
Let $x_j = c_1(L_j)$, $\sigma_j = \sigma_j(x_1, \cdots, x_n)$
the $j$-th elementary symmetric polynomial in $x_1, \cdots, x_n$.
Then $(\ref{eqn:flag})$ is equivalent to 
$\sigma_j = 0$, $j = 1, \cdots, n$.
Now we set $x_{j, h} = c_1(L_j)_h$, $\sigma_{j, h} =
\sigma_j(x_{j, h}, \cdots, x_{n, h})_h$, the $j$-th elementary
symmetric polynomial computed by quantum multiplications.
Then from quantum Chern-Weil theory, we have
$\sigma_{j, h} = 0$, $j = 1, \cdots, n$.
As in the projective space case, this does not help much 
in writing down the quantum multiplications in quantum de Rham
cohomology.

We identify $F_{n+1}$ with $X = U(n)/T^n$, where $T^n$ is
the diagonal subgroup.
Let $x_0$ denote the point in $X$ which corresponds to $T^n$.
Then we can identify $T_{x_0}X$ with
$$\fn = \{ (a_{ij}) \in M(n \times n, {\Bbb C}): a_{jj} = 0,
j = 1, \cdots, n, a_{ij} = - \bar{a}_{ji} \}.$$ 
Since $H^*_{dR}(U(n)/T^n) \cong H^*(\Omega(U(n)/T^n)^{U(n)}, d)$,
every de Rham class can be represented by a $U(n)$-invariant closed
form.
In particular, by choosing a $U(n)$-invariant connection, 
the first Chern class of the line bundle $L_j$,  
can be represented by a $U(n)$-invariant closed form $\alpha_j$, 
for each $j = 1, \cdots n$.
The calculations of quantum de Rham cohomology for 
$U(n)$-invariant symplectic structures on $U(n)/T^n$
will be indicated later.
\end{example}

\begin{example} (Generalized flag manifolds)
We identify $F_{n+1}$ with $X = U(n)/T^n$, where $T^n$ is
the diagonal subgroup.
This reveals the fact that complex flag manifold is 
a special example of an important class of 
Fano manifolds used in Borel-Weil-Bott theory. 
In general, let $G^{\Bbb C}$ be a semisimple Lie group 
over ${\Bbb C}$,
and $B$ be a Borel subgroup of $G^{\Bbb C}$, 
then $G^{\Bbb C}/B$ is a projective variety (see Borel \cite{Bor2}). 
In fact, let $G$ be the maximal compact subgroup of $G^{\Bbb C}$, 
and $T$ the maximal torus of $G$, then $G/B \cong G/T$. 
Use a  $G$-invariant inner product, e.g. 
the negative of Killing form on $\fg$,
one gets a decomposition
$$\fg = \ft \oplus \fn.$$
one can identify $T_{x_0}(G/T)$ with $\fn$, 
where $x_0$ is the point in $X = G/T$ corresponding to $T$.
It is straightforward to find compatible almost complex structure $J$
and symplectic structure $\omega$ on $T_{x_0}X$ which are
invariant under the action of the Weyl group $W$.
Use the translation by action of $G$ to extend $J$ and $\omega$
to $G/T$, we then get a homogeneous K\"{a}hler manifold,
which is Fano.
Each weight $\lambda \in \ft^*$ of $G$ determines 
a representation $\phi_{\lambda}: T \rightarrow S^1$. 
Since $\pi: G \rightarrow G/T$ is a principal $T$-bundle,
we then get 
an associated line bundle $L_{\lambda}$ from $\phi_{\lambda}$. 
Let $\lambda_1, \cdots, \lambda_k$ be the set of simple roots of
$G$, and $x_j = c_1(L_{\lambda_j})$. 
If $p \in I(G) = S({\frak k}^*)^G$ is an invariant polynomial, 
then by the isomorphism $S({\frak G}^*)^G \cong S({\frak t}^*)^W$,
we can identify $p$ with a polynomial on $\frak t$,
which is invariant under the action of the Weyl group.
It turn out that $p(x_1, \cdots, x_k) = 0$, if $p(0) = 0$. 
Let $I_+(G)$ be the ideal generated by $p \in S(\ft^*)^T$ such that
$p(0) = 0$. 
Borel \cite{Bor1} used the degeneracy of the Leray spectral sequence 
at $E_2$ of the fibration $G/T \rightarrow BG \rightarrow BT$ 
to show that the cohomology of $X$ is isomorphic to 
$$S({\frak t}^*)/I_+(G).$$ 
When $G = U(n)$, $T = T^n$ the diagonal subgroup,  
$G/T$ is diffeomorphic to the complex flag manifold,
one recovers the result of the last example.

To compute the quantum de Rham cohomology, 
we  use the fact $H^*_{dR}(G/T) = H^*(\Omega(G/T)^G, d)$,
which implies that de Rham classes of $G/T$
can be represented by  $G$-invariant forms on $G/T$.
Such forms are determined by their values at $x_0$.
Since $T_{x_0}(G/T) \cong\fn$,
$\Omega(G/T)^G \cong \Lambda(\fn^*)^T$.
For any element $\xi \in \fg$, let $X_{\xi}$ be the fundamental
vector field of $\xi$ on $X = G/T$. 
For any weight $\lambda$, define a $G$-invariant 
$2$-form  $\omega_{\lambda}$ by setting
$$\omega_{\lambda}(X_{\xi}, X_{\eta}) = \lambda([\xi, \eta]),$$ 
at $x_0$, for two fundamental vector fields 
$X_{\xi}$ and $X_{\eta}$.
Following  Lemma 8.67  and Lemma 8.68 in Besse \cite{Bes},
it is easy to see that $\omega_{\lambda}$ is closed.
Let $\lambda_1, \cdots, \lambda_r$ be the simple roots of $\fg$,
then they form a basis of $\ft^*$. 
Borel's result indicates that $\omega_{\lambda_1}, \cdots, 
\omega_{\lambda_r}$ generates $H^*_{dR}(G/T)$.
We will consider $G$-invariant symplectic forms 
which will  be given below.
\end{example}

\begin{example} (Coadjoint orbits) 
The coadjoint orbits of a compact connect Lie group $G$ 
are parameterized by $\ft^*/W$, or equivalently, 
a closed Weyl chamber $\overline{C}$.
For any $\lambda \in \overline{C}$,
let ${\cal O}_{\lambda}$ denote the orbit of $\lambda$.
For $\lambda$ in the interior of $\overline{C}$, 
${\cal O}_{\lambda}$ is diffeomorphic to $G/T$;
for $\lambda$ on the wall of the Weyl chamber, 
there is a fibration of $G/T$ over ${\cal  O}_{\cal O}$
(Besse \cite{Bes}, Proposition 8.116). 
Besse \cite{Bes}, $\S 8H$ shows complex flag
manifolds, partial flag manifolds, Grassmannians, etc., 
are all examples of coadjoint orbits.
It can be shown that $\omega_{\lambda}$ is 
a  $G$-invariant symplectic forms on ${\cal O}_{\lambda}$.
It is called Kirillov-Kostant-Souriau form.
It is straightforward   to explicitly write down the 
quantum exterior multiplication on 
$H^*_{dR}({\cal O}_{\lambda}, \omega_{\lambda})$,
since we can do it on the tangent space of one point.
\end{example}

\end{document}